\documentclass[twocolumn,10pt]{asme2ej}

\usepackage{graphicx}

\usepackage{amsmath}
\usepackage{amssymb}
\usepackage{mathrsfs}

\usepackage{todonotes}

\usepackage{siunitx}
\newcommand{\dimSpace}{d}
\newcommand{\dimSol}{k}
\newcommand{\numParams}{{n_p}}
\newcommand{\numGeoPar}{{n_g}}
\newcommand{\sol}{\mathbf{u}}
\newcommand{\state}{\mathbf{y}}
\newcommand{\pt}{\mathbf{x}}
\newcommand{\parVec}{\boldsymbol{\mu}_p}
\newcommand{\geoVec}{\boldsymbol{\mu}_g}
\newcommand{\geoSpace}{\mathcal{G}}
\newcommand{\parSpace}{\mathcal{P}}
\newcommand{\geoMapping}{{P_g}}

\newcommand{\diffOperator}{\mathcal{L}}
\newcommand{\boundOperator}{\mathcal{B}}
\newcommand{\obsOperator}{\mathcal{U}}

\newcommand{\numTrain}{N_{sn}}
\newcommand{\numPoints}[1]{N_{pt}^{#1}}
\newcommand{\ptObs}[2]{\pt^{#1}_{#2}}

\newcommand{\NN}{\mathcal{NN}}
\newcommand{\dist}{d}

\newcommand{\annParams}{\mathbf{w}}
\newcommand{\Loss}{\mathcal{J}}
\newcommand{\RegTerm}{\mathcal{R}}

\newcommand{\UCsys}{\Phi_\geoSpace}
\newcommand{\DomainRef}{\widehat{\Omega}}
\newcommand{\ptRef}{\widehat{\mathbf{x}}}

\newcommand{\Rey}{\mathbb{R}\mathrm{e}}
\newcommand{\pot}{\psi}
\newcommand{\vel}{\mathbf{v}}
\newcommand{\prs}{p}

\newcommand{\UCcoroLR}{\psi_{\mathrm{LR}}}
\newcommand{\UCcoroTD}{\psi_{\mathrm{TD}}}

\title{Universal Solution Manifold Networks (USM-Nets): non-intrusive mesh-free surrogate models for problems in variable domains}

\author{Francesco Regazzoni 
    \affiliation{
        MOX - Department of Mathematics, \\Politecnico di Milano, Milan, Italy. \\
        \texttt{francesco.regazzoni@polimi.it}
    }	
}

\author{Stefano Pagani
    \affiliation{
        MOX - Department of Mathematics, \\Politecnico di Milano, Milan, Italy. \\
        \texttt{stefano.pagani@polimi.it}
    }	
}

\author{Alfio Quarteroni
    \affiliation{
        MOX-Department of Mathematics, Politecnico di Milano, Milan, Italy. \\
        Institute of Mathematics, EPFL, Lausanne, Switzerland (\textit{Professor Emeritus}) \\
        \texttt{alfio.quarteroni@polimi.it}
    }	
}

\begin{document}

\maketitle    

\begin{abstract}
{\it 
We introduce Universal Solution Manifold Network (USM-Net), a novel surrogate model, based on Artificial Neural Networks (ANNs), which applies to differential problems whose solution depends on physical and geometrical parameters. Our method employs a mesh-less architecture, thus overcoming the limitations associated with image segmentation and mesh generation required by traditional discretization methods. Indeed, we encode geometrical variability through scalar landmarks, such as coordinates of points of interest. In biomedical applications, these landmarks can be inexpensively processed from clinical images. Our approach is non-intrusive and modular, as we select a data-driven loss function. The latter can also be modified by considering additional constraints, thus leveraging available physical knowledge. Our approach can also accommodate a universal coordinate system, which supports the USM-Net in learning the correspondence between points belonging to different geometries, boosting prediction accuracy on unobserved geometries. Finally, we present two numerical test cases in computational fluid dynamics involving variable Reynolds numbers as well as computational domains of variable shape. The results show that our method allows for inexpensive but accurate approximations of velocity and pressure, avoiding computationally expensive image segmentation, mesh generation, or re-training for every new instance of physical parameters and shape of the domain.
}
\end{abstract}

\section{Introduction} \label{sec:intro}

Models and methods in scientific computing and machine learning enable the extraction of relevant knowledge from available data \cite{alber2019integrating}.
Data can represent, e.g., physical coefficients, geometrical factors, boundary or initial conditions.
This is the case of computational fluid dynamics, a remarkable instance, especially when addressing problems in aerodynamics, such as the design of vehicles, and in biomedical engineering, such as patient-specific analysis of blood flow dynamics \cite{anderson1995computational,parolini2005mathematical,formaggia2010cardiovascular,brunton2020machine}.

\begin{figure*}
    \centering
    \includegraphics[width=\textwidth]{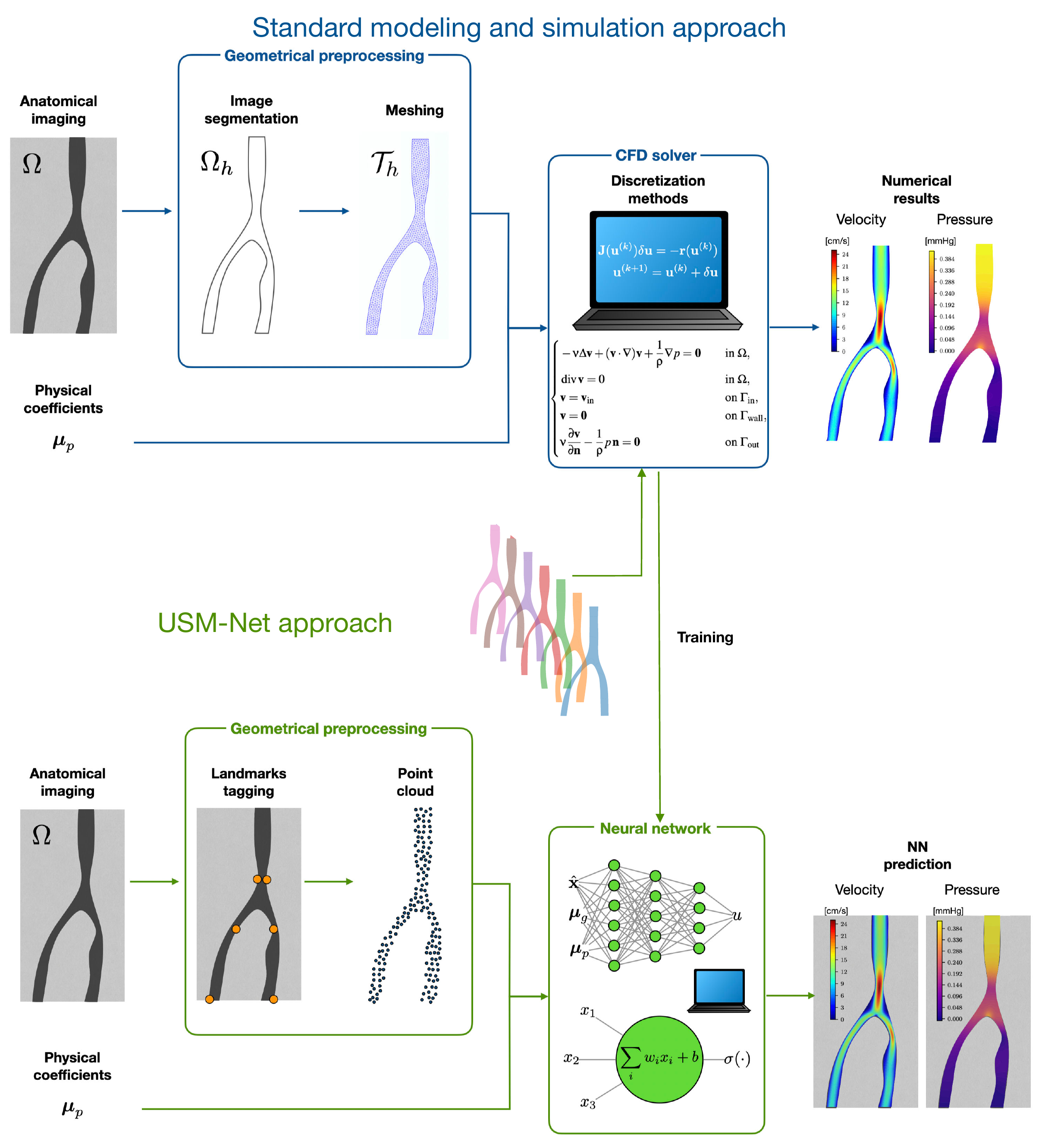}\\
    \caption{Comparison between the standard modeling and simulation approach (top) and the USM-Nets approach (bottom) for a use case of clinical interest, namely the prediction of blood flow and pressure within a coronary bifurcation.
    The former approach requires a geometric preprocessing phase, which consists in segmenting the patient's clinical images to extract a computational domain $\Omega_h$. A partition of the latter into a set of cells (in the figure, into triangles) constitutes the computational mesh $\mathcal{T}_h$.
    The numerical solution is then obtained by the discretization of the Navier-Stokes equations on the computational mesh $\mathcal{T}_h$, solved through suitable computer-based algorithms.
    Our proposed USM-Net approach lightens both the geometric preprocessing and the solution approximation steps. The first one consists solely of landmarks extraction and (mesh-less) cloud of query points generation. Finally, the solution is obtained by evaluating the neural network as many times as the number of query points.
    The USM-Net can be trained either from experimentally measured data or from synthetic data generated using the Navier-Stokes solver.
    }
    \label{fig:approach_comparison}
\end{figure*}

The standard approach to modeling and simulation (Fig.~\ref{fig:approach_comparison}, top) requires the construction of a computational mesh, a partition of the given geometry in simple elements, e.g. tetrahedral or hexahedral cells.
In biomedical applications, the construction of the computational mesh requires a preliminary step of segmentation, i.e. the extraction of the boundaries of the computational domain from medical images, such as those derived from magnetic resonance or computerized tomography.
Then, discretization methods (like Finite Elements and Finite Difference methods \cite{quarteroni2017numerical}) assemble on the elements of the mesh a suitable approximation of the operators associated with the partial differential equations (PDEs).
Unfortunately, changes in shape require the re-execution of the entire process, necessitating the reallocation of significant computational resources.

For this reason, some computational approaches, like isogeometric analysis (IGA) and shape models, are designed to avoid the regeneration of a new computational mesh when a change of geometry occurs.
IGA achieves this thanks to the use of non-uniform rational B-splines (NURBS) that exactly match CAD geometries, usually adopted in an industrial context~\cite{hughes2005isogeometric}.
Similarly, geometrical shape models describe the possible variations of geometry through a limited number of parameters.
Among the most diffused shape models we recall free-form deformation (FFD)~\cite{sederberg1986free,lamousin1994nurbs}, radial basis function (RBF)~\cite{buhmann2000radial}, and statistical shape models \cite{heimann2009statistical}, based e.g. on principal component analysis (PCA) \cite{jolliffe2016principal}.
Firstly conceived in computer graphics, FFD is a technique that surrounds the object with a lattice of control points.
Their motion drives the deformation of all object points through a polynomial interpolation.
RBFs define a parametrized map describing the deformations of the geometry from a small number of selected control points as well.
Compared to the FFD approach, the control points position is not constrained on a lattice but user-selected from application to application.
Finally, PCA-based models describe a collection of similar geometries through weighted modes describing the principal geometrical variations with respect to a mean shape.
These approaches enable the construction of reduced-order models (ROMs) for problems with parameterized geometry, which provide a computationally efficient approximation of the solution for many different choices of the geometrical parameters.
In this context, empirical interpolation methods build efficient approximations of differential operators avoiding computationally costly reassembly.
Among the main applications of parametrized ROMs \cite{quarteroni2011certified,carlberg2011efficient,quarteroni2015reduced,hesthaven2016certified,peherstorfer2016data,taira2017modal}, we mention those related to computational fluid dynamics using FFD \cite{samareh2004aerodynamic,lassila2010parametric}, RBF \cite{morris2008cfd,rendall2008unified,manzoni2012model} and PCA \cite{sangalli2009case}.

Besides introducing a geometrical error, the critical aspect of shape models is that, when deforming a fixed mesh, the elements may encounter such deformations that the discretized problem becomes ill-conditioned. This might considerably limit the ability of the method to explore the geometric variability of the problem.
To overcome the limitations of mesh-based approaches, mesh-free and particle methods construct a discretization of the geometry formed solely by a collection of points, relaxing the constraints given by the connectivity \cite{li2002meshfree}. These methods bring numerous advantages in terms of managing geometrical accuracy (even with discontinuities), imposing large deformations, adaptively refining, and code parallelization. However, they generate full-order models (FOMs) with limited computational efficiency.

In this paper, we introduce a new surrogate modeling technique based on artificial neural networks (ANNs), that can learn the solution manifold of a given parametrized PDE.
Notably, these surrogate models are not tailored on a given domain, but they account for the influence of geometry on the solution.
We leverage the capabilities of approximating arbitrary complex functions with inexpensive output evaluation provided by ANNs, which are indeed \textit{universal} approximators of several families of functions, including continuous functions \cite{cybenko1989approximation} and Sobolev spaces \cite{mhaskar_neural_1996,montanelli2019new}.
Based on these results, the surrogate models we propose are, in principle, able to approximate the solution manifold of a given differential problem with arbitrary accuracy and universally with respect to the variation of domain geometry and physical parameters.
For this reason, we name them \textit{Universal Solution Manifold Networks} (USM-Nets).
We train these networks using subsamples of solutions snapshots obtained by varying the physical parameters and the domain geometry (either synthetically generated through a FOM or experimentally collected). Therefore, the accuracy of the predictions of a specific USM-Net would directly depend on the richness of the training set used. However, the architecture we propose has the potential to accurately approximate the solution universally with respect to domain and physical parameters, which is not the case for methods that are constrained to a predetermined parameterization of the solution manifold.

The design of a USM-Net avoids complex geometrical preprocessing (comprising segmentation and construction of the computational mesh; see Fig.~\ref{fig:approach_comparison}).
We encode geometrical variability in a few \textit{geometrical landmarks}, a finite number of scalar indicators, inexpensively processed from the input image.
In the most basic case, landmarks are the coordinates of specific reference points where we impose a correspondence between geometries.
Landmarks might identify inlets, bifurcations, or other specific structures, and they provide, together with the physical parameters, the input to the ANN.

%

Our approach is modular, thanks to the possibility of arbitrarily defining the loss function. Besides penalizing the misfit with the available data, during the training of the network we can enforce assumptions of regularity (imposing, e.g., weights penalization), the initial or boundary conditions, or the fulfillment of an equation in a strong (differential) form. The latter would represent a generalization of the so-called physics-informed neural network \cite{raissi2019physics,raissi2020hidden,regazzoni2021physics}, avoiding new executions of the training process for each geometry variation.

The outline of this paper is as follows.
First, in Sec.~\ref{sec:methods} we present the notation used throughout this paper and the proposed methods.
Then, in Sec.~\ref{sec:test-cases} we present two test cases and we describe how our proposed methods can be applied to them.
In Sec.~\ref{sec:results} we present the numerical results and in Sec.~\ref{sec:conclusions}, finally, we draw our conclusions.

\section{Methods} \label{sec:methods}

In this section, we first introduce the notation used throughout this paper.
Then, we present our proposed USM-Net method.

\subsection{Problem setting}

We consider a space-dependent physical quantity $\sol(\pt)$, defined on a domain $\Omega \subset \mathbb{R}^\dimSpace$, where typically $d=2, 3$.
For the sake of generality, we denote $\sol(\pt) \in \mathbb{R}^{\dimSol}$ as a vector field even though in certain cases it could be a scalar field (that is $\dimSol=1$).
For example, $\sol(\pt)$ may correspond to the blood velocity field in a vessel, the displacement field of soft tissue, or the pressure field of air around a body.
Very often, the quantity $\sol(\pt)$ depends on a set of physical parameters $\parVec \in \parSpace$, where $\parSpace \subset \mathbb{R}^{\numParams}$ is the parameter space.
The parameters $\parVec$ characterize the physical processes that determine $\sol$.
For example, the velocity of a fluid depends on its viscosity, while the displacement of biological tissue depends on its stiffness moduli and the applied force.
Furthermore, the field $\sol(\pt)$ depends on the domain $\Omega$ itself.
In many practical applications, we are interested in the family of fields $\sol(\pt)$ obtained by varying the domain $\Omega$ in a given set, denoted by $\geoSpace$:
\begin{equation*}
    \geoSpace \subset \{ \Omega \subset \mathbb{R}^\dimSpace, \text{ open and bounded}\}.
\end{equation*}
To stress the dependence of $\sol(\pt)$ on both $\parVec$ and $\Omega$, we will henceforth write $\sol(\pt; \parVec, \Omega)$.

Often, the physical processes that determine $\sol(\pt; \parVec, \Omega)$ can be described by a mathematical model assuming the form of a boundary-value problem, that is
\begin{equation} \label{eqn:PDE_generic}
    \left\{
    \begin{aligned}
        \diffOperator(\state, \parVec) = \mathbf{0}  \qquad& \text{for $\pt \in \Omega$}, \\
        \boundOperator(\state, \parVec) = \mathbf{0} \qquad& \text{for $\pt \in \partial\Omega$}, \\
        \sol = \obsOperator(\state, \parVec)  \qquad& \text{for $\pt \in \Omega$}, \\
    \end{aligned}
    \right.
\end{equation}
where $\state$ in the state variable; $\diffOperator$ and $\boundOperator$ are the operators associated with the differential equation and with the boundary conditions, respectively; $\obsOperator$ is the observation operator.
Since the state $\state$ is instrumental for obtaining $\sol$, in what follows we will refer to $\sol$ and not to $\state$ as the \textit{solution} of the FOM.
If the differential problem \eqref{eqn:PDE_generic} is well-posed, then, given $\parVec \in \parSpace$ and $\Omega \in \geoSpace$ there exists a unique solution $\sol(\cdot; \parVec, \Omega) \colon \Omega \to \mathbb{R}^{\dimSol}$.
This solution can be numerically approximated through a FOM based e.g. on Finite Differences or Finite Elements \cite{quarteroni2017numerical}.

The goal of this paper is to build an emulator that surrogates the solution map $(\parVec, \Omega) \mapsto \sol(\cdot; \parVec, \Omega)$, that is a model that, for any given value of the geometrical parameters $\parVec$ and any given geometry $\Omega$, provides an approximation of the corresponding solution $\sol$.

\subsection{USM-Net}

A USM-Net is an ANN-based model, trained on a data-driven basis, that surrogates the solution map $(\parVec, \Omega) \mapsto \sol(\cdot; \parVec, \Omega)$, without using any FOM.
In the most common scenario, a FOM that describes the physical process is available, and it is used to generate the data needed to train the USM-Net.
Still, the training of a USM-Net is also possible for physical problems in which the FOM is unknown, provided that sufficient experimental measurements are available.
We present two versions of USM-Net:
\begin{enumerate}
    \item {PC-USM-Net}, when the solution is represented in terms of the \textit{physical coordinates}, that is $\pt \in \Omega$ (see Sec.~\ref{sec:methods:PC-USM-Nets});
    \item {UC-USM-Net}, with the solution represented now by passing through of a system of \textit{universal coordinates}, that will be defined in Sec.~\ref{sec:methods:UC}.
\end{enumerate}

\subsubsection{PC-USM-Net} \label{sec:methods:PC-USM-Nets}

\begin{figure}[ht]
    \includegraphics[width=\columnwidth]{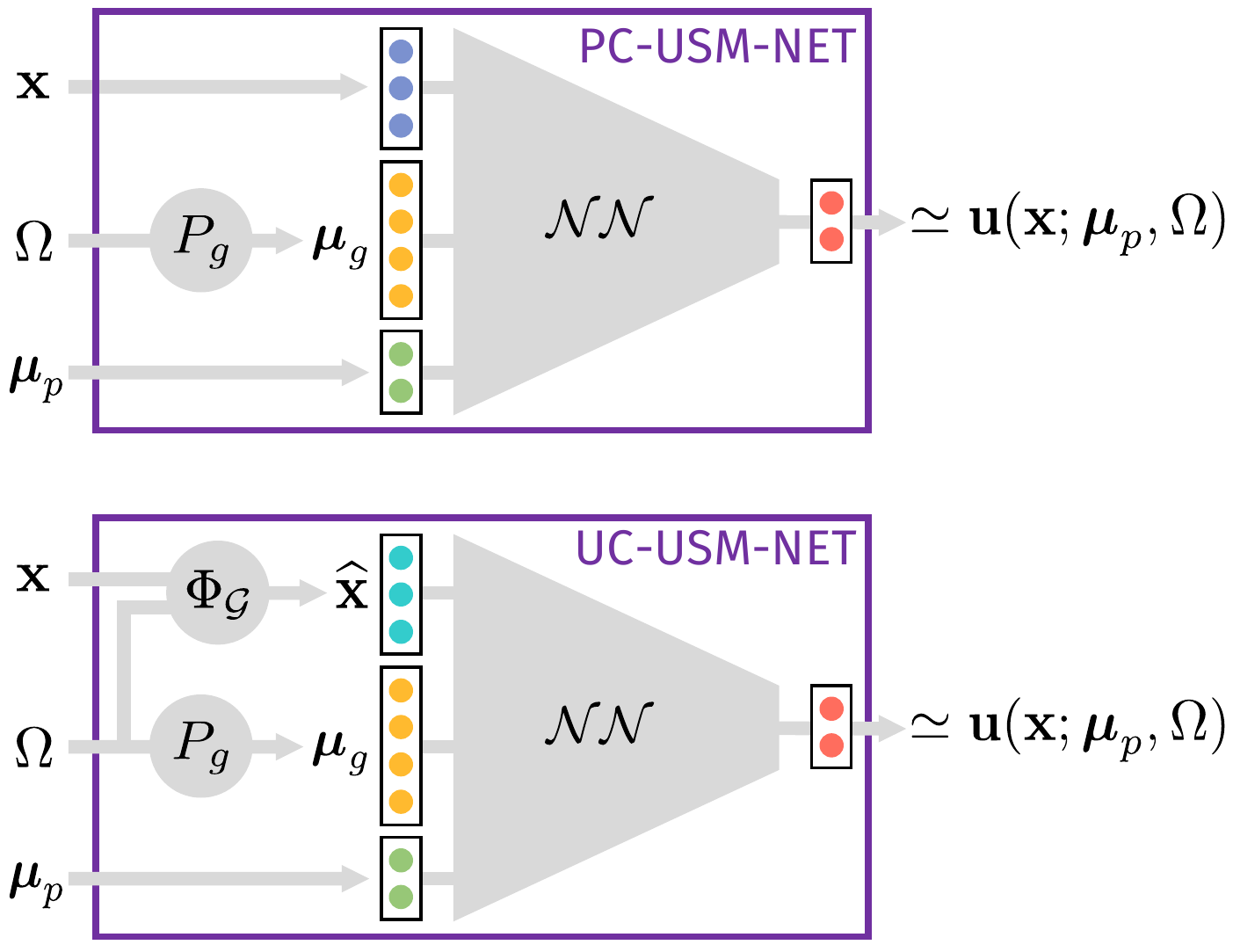}
    \caption{Architecture of a PC-USM-Net (top) and of a UC-USM-Net (bottom).}
    \label{fig:architecture_merged}
\end{figure}

The PC-USM-Net architecture is represented in Fig.~\ref{fig:architecture_merged} (top). It consists of an ANN, typically a fully connected ANN (FCNN), whose input is obtained by stacking three vectors:
\begin{enumerate}
    \item the query point $\pt$, that is the coordinate of the point where the solution $\sol(\pt)$ is sought;
    \item the physical parameters $\parVec$;
    \item a set of geometrical landmarks associated with the domain at hand, that is $\geoVec = \geoMapping(\Omega)$, that typically represent the coordinates of key points of the domain (such as inlets, bifurcation points, etc.) or geometrical measures (such as diameters, thicknesses, etc.).
    In Sec.~\ref{sec:methods:landmarks} we will elaborate on possible choices for the function $\geoMapping$.
\end{enumerate}
The output of the FCNN is an approximation of the solution $\sol$ at the query point $\pt$.
More precisely, denoting by $\NN$ the FCNN and by $\annParams$ its trainable parameters (weights and biases), we have:
\begin{equation*}
    \sol(\pt; \parVec, \Omega) \simeq \NN(\pt, \parVec, \geoMapping(\Omega); \annParams).
\end{equation*}
Hence, $\NN$ features $\dimSpace + \numParams + \numGeoPar$ input neurons and $\dimSol$ output neurons.

\subsubsection{UC-USM-Net} \label{sec:methods:UC}

As anticipated, a PC-USM-Net has a universality character with respect to domains, i.e. a single network is used to represent solutions in different geometries.
However, a given point with physical coordinate $\pt$, given as input to the PC-USM-Net, can play a different \textit{role} for different geometries.
For example, a point $\pt$ that for one geometry $\Omega^1$ belongs to the boundary of the domain, for another geometry $\Omega^2$ could be internal to the domain, and for yet another $\Omega^3$ could even be external.
Therefore, we propose an evolution of PC-USM-Net aimed at capturing more effectively the correspondence between points among geometries.

To achieve this goal, we rely on a \textit{universal coordinates (UC) system} for $\geoSpace$.
A UC system is a map $\UCsys$ that, to any domain $\Omega$ and to any point $\pt \in \Omega$, associates a point $\ptRef \in \DomainRef$ belonging to a \textit{reference domain} $\DomainRef \subset \mathbb{R}^{\dimSpace}$.
More precisely, the reference coordinate is obtained as $\ptRef = \UCsys(\pt,\Omega)$.
We require that the application $\UCsys(\cdot, \Omega) \colon \Omega \to \DomainRef$ be a continuous bijection for any $\Omega \in \geoSpace$.
Hence, the UC system $\UCsys$ defines a coordinate transformation that maps each domain $\Omega \in \geoSpace$ into the reference one $\DomainRef$.
In Sec.~\ref{sec:test-cases} we will present two concrete examples of UC systems.

Whenever a UC system is available, it can augment a PC-USM-Net.
The enhanced version of a PC-USM-Net, called UC-USM-Net, is obtained by giving as input to $\NN$ the reference coordinate $\ptRef = \UCsys(\pt, \Omega)$ instead of the physical one $\pt \in \Omega$.
More precisely, the surrogate model is defined as
\begin{equation} \label{eqn:UC-USM-Net}
    \sol(\pt; \parVec, \Omega) \simeq \NN(\UCsys(\pt, \Omega), \parVec, \geoMapping(\Omega); \annParams).
\end{equation}
The resulting architecture is represented in Fig.~\ref{fig:architecture_merged} (bottom).


We remark that UC-USM-Net is a generalization of PC-USM-Net, as the latter can be obtained from the former by setting $\UCsys$ equal to the identity function, that is by setting $\ptRef = \pt$.
For this reason, from now on we will consider without loss of generality the surrogate model of Eq.~\eqref{eqn:UC-USM-Net}.

As we will show in the results section, UC-USM-Nets allow to improve the generalization accuracy of PC-USM-Nets, that is the accuracy of predictions for physical parameters and geometries not included in the training set, by providing geometrical prior knowledge during training.
Moreover, we will show that, besides helping the ANN to link together points of different geometries, a UC system might provide details of the geometry that are not captured by the landmarks.

We remark that, in practical applications, both the injectivity and the surjectivity requirements of the map $\UCsys(\cdot, \Omega) \colon \Omega \to \DomainRef$ can be relaxed.

\subsection{Geometrical landmarks} \label{sec:methods:landmarks}

In order to build a surrogate model that learns an approximation of the solution map, we introduce a low-dimensional description of the geometry.
In particular, we construct an operator $\geoMapping \colon \geoSpace \mapsto \mathbb{R}^{\numGeoPar}$ that, to any given computational domain $\Omega \in \geoSpace$, associates a finite number (say $\numGeoPar$) of geometrical landmarks $\geoVec = \geoMapping(\Omega) \in \mathbb{R}^{\numGeoPar}$, which are suitable to provide a compact description of $\Omega$.
Depending on the structure of $\geoSpace$, different strategies can be followed to define the operator $\geoMapping$.
\begin{enumerate}
    \item In case an explicit parametrization of the elements of the space $\geoSpace$ is available, we define $\geoMapping$ in such a way that the landmarks $\geoVec = \geoMapping(\Omega)$ coincide with the geometrical parameters themselves.
    An example is provided in Test Case~1.

    \item If such a parameterization is not available (as it is in many cases of practical interest), a straightforward choice is to take the coordinates of key points in the domain as landmarks.
    An example is provided in Test Case~2.

    \item Other more sophisticated techniques can be used to obtain a low-dimensional description of the computational domains.
    For example, the geometrical landmarks can be defined as the first, more relevant, coefficients associated with a POD analysis of a finite subset of $\geoSpace$ (shape model).
    Entering into the details of this or other techniques is beyond the scope of this paper.
    The method proposed in this paper is indeed general, as it is built on top of the different techniques that can be used to define the map.

\end{enumerate}
In general, our method does not require the operator $\geoMapping$ to be invertible.
As a matter of fact, $\geoMapping$ is invertible only when an explicit parametrization of the space $\geoSpace$ is available.
In fact, we allow for the case when two different geometries $\Omega^1 \neq \Omega^2$, both belonging to $\geoSpace$, are associated to identical landmarks (i.e. $\geoMapping(\Omega^1) = \geoMapping(\Omega^2)$).
Since the geometrical landmarks characterize the variability of the geometry, a good design of $\geoMapping$ requires that the condition $\geoMapping(\Omega^1) = \geoMapping(\Omega^2)$ implies that $\Omega^1$ and $\Omega^2$ are {\textit{minimally} different}.

\subsection{Training a USM-Net}

To train a USM-Net (that is, either a PC-USM-Net or a UC-USM-Net), we require the output of problem~\eqref{eqn:PDE_generic} for several pairs $(\parVec, \Omega) \in \parSpace \times \geoSpace$.
These solutions, called \textit{snapshots}, are typically obtained through the FOM, a high-fidelity numerical solver of Eq.~\eqref{eqn:PDE_generic}, based e.g. on Finite Elements of Finite Differences \cite{quarteroni2017numerical}.
Yet, as our method is non-intrusive and does not require any knowledge of equation~\eqref{eqn:PDE_generic}, snapshots can also be derived from a black-box solver, or even from experimental measurements.

We consider a collection of $\numTrain$ snapshots, associated with $\parVec^{i} \in \parSpace$ and $\Omega^{i} \in \geoSpace$, for $i = 1, \dots, \numTrain$.
For any snapshot, then, we consider a number of observations in a set of points belonging to $\Omega^{i}$.
In case of high variability of the geometries in $\geoSpace$, the resolution of the FOM typically requires the generation of different meshes, without a one-to-one correspondence of the nodes.
Therefore, to guarantee generality, we consider the case where each snapshot has a potentially different number of observation points.
Specifically, we denote by
$\{\ptObs{i}{j}, \, j=1,\dots,\numPoints{i}\}$
the set of observation points associated with the $i$-th snapshot.
In conclusion, the training dataset consists of the following set
\begin{equation*}
    \{
        \parVec^{i},
        \Omega^{i},
        \{
            \sol(\ptObs{i}{j}; \parVec^{i}, \Omega^{i})
        \}_{j=1}^{\numPoints{i}}
    \}_{i=1}^{\numTrain}.
\end{equation*}
Training the USM-Net requires solving the following minimization problem
\begin{equation*}
    \widehat{\annParams} = \underset{\annParams}{\operatorname{argmin}} \; \Loss(\annParams).
\end{equation*}
The loss function $\Loss$ is given by the misfit between snapshot data and predictions, plus (optionally) a regularization term $\RegTerm$:
\begin{equation}\label{eqn:loss}
    \Loss(\annParams) =
            \frac{1}{\numTrain}\sum_{i=1}^{\numTrain} \frac{1}{ \numPoints{i}} \sum_{j=1}^{\numPoints{i}}
            \dist( \sol_j^i, \widetilde{\sol}_j^i)+ \RegTerm(\annParams),
\end{equation}
having defined
\begin{equation*}
    \begin{split}
        \sol_j^i &= \sol(\ptObs{i}{j}; \parVec^{i}, \Omega^{i}), \\
        \widetilde{\sol}_j^i &=
        \NN(\UCsys(\ptObs{i}{j}, \Omega^{i}), \parVec^{i}, \geoMapping(\Omega^{i}); \annParams)
    \end{split}
\end{equation*}
and where $\dist(\cdot, \cdot)$ is a suitable discrepancy metric (typically, $\dist(\sol, \tilde{\sol}) = \| \sol -\tilde{\sol}\|^2$).
Standard techniques, such as Tikhonov or LASSO regularization, can be used for the regularization term $\RegTerm$.
Additionally, $\RegTerm$ can be augmented by suitable terms informing the USM-Net of physical knowledge available on the solution (see also Sec.~\ref{sec:methods:gray-box}).
An example in this sense will be shown in Sec.~\ref{sec:results}.

\subsection{Grey-box USM-Net} \label{sec:methods:gray-box}

So far, we have presented USM-Nets as fully non-intrusive (black-box) surrogate modeling techniques.
Still, physical knowledge can be optionally embedded into their construction.
Indeed, the training process can be augmented by informing the network either of physical constraints (such as conservation principles, symmetry properties, or the positivity of the solution) or of the differential equations and boundary conditions that characterize the problem.
We distinguish between \textit{weak imposition} and \textit{strong imposition} of the physical knowledge.

\paragraph{Weak imposition.} Prior knowledge on the solution map can be enforced through the regularization term $\RegTerm$ of the loss function of \eqref{eqn:loss}.
$\RegTerm$ can include the norm of the residual of the FOM equations and boundary conditions evaluated in a collection of collocation points, as done in the training of Physics Informed Neural Networks \cite{raissi2019physics}.
Similarly, other physical constraints can be rephrased in terms of minimization of a regularization term $\RegTerm$.
Thanks to Automatic Differentiation, the inclusion of differential operators in the term $\RegTerm$ does not introduce severe implementation efforts, even if it will slow down the training process.
Therefore, the user should wisely balance the advantages and disadvantages of introducing such a term.
An example of the weak imposition of the boundary conditions is presented in Test Case 1 (Sec.~\ref{sec:test-cases:cavity}).

\paragraph{Strong imposition.} Alternatively, we can enforce physical constraints by defining an architecture $\NN$ that satisfies them by construction.
We now give a brief list of examples:
\begin{enumerate}
    \item Non-negativity of the solution can be enforced by introducing after the FCNN a further layer that applies an operator with nonnegative output, such as $(\cdot)^2$ or $|\cdot|$.
    In other terms, we perform a composition between the FCNN and the nonnegative operator.
    \item Symmetry w.r.t. a given input coordinate can be enforced, e.g., by introducing an input layer that pre-processes the corresponding input through an even function, such as $|\cdot|$.
    As in the previous point, this corresponds to performing a composition between the even function and the FCNN.
    \item Dirichlet boundary conditions on a portion of the boundary (i.e. $\sol(\pt) = \sol_D$ on $\Gamma_{D} \subset \partial\Omega$) can be strongly enforced by introducing a multiplicative layer after the FCNN that multiplies the solution by a mask function $\Phi_{\text{BC}}(\pt, \Omega)$, such that $\Phi_{\text{BC}} = 0$ on $\Gamma_{D}$ and $\Phi_{\text{BC}} \neq 0$ elsewhere, and sums the datum $\sol_D$.
    \item Solenoidality of the solution ($\nabla \cdot {\sol} = 0$), a common requirement in fluid dynamics, can be enforced by interpreting the FCNN output as the flow field potential and introducing an output layer that returns its curl.
    An example of this technique is described in Sec.~\ref{sec:test-cases:cavity}.
\end{enumerate}

\subsection{Notes about implementation}

From the implementation point of view, a few cautions are needed to make the training of USM-Nets computationally light.
First of all, despite the application of the $\geoMapping$ map in Fig.~\ref{fig:architecture_merged} is indicated as being an integral part of the USM-Net, the landmarks $\geoVec^i = \geoMapping(\Omega^i)$ can be pre-calculated before the training.
Similarly, in the case a UC system is employed, the coordinate transformation $\ptRef_j^i = \UCsys(\pt_j^i, \Omega^i)$ for any $i$ and $j$ can be performed offline, at a stage prior to training.
In this manner, we set up an augmented dataset consisting of:
\begin{equation*}
    \begin{aligned}
    (&\ptRef_1^1,&& \geoVec^1, && \parVec^1), \qquad && \sol_1^1 \\
    (&\ptRef_2^1,&& \geoVec^1, && \parVec^1), \qquad && \sol_2^1 \\
    &\vdots     && \vdots     && \vdots     && \vdots   \\
    (&\ptRef_{\numPoints{1}}^1,&& \geoVec^1 && \parVec^1), \qquad && \sol_{\numPoints{1}}^1 \\
    (&\ptRef_1^2,&& \geoVec^2, && \parVec^2), \qquad && \sol_1^2 \\
    &\vdots     && \vdots     && \vdots     && \vdots   \\
\end{aligned}
\end{equation*}
and $\NN$ is trained to fit the map from the first three columns to the last one.

Once the $\NN$ is trained, it can be used to approximate the solution for unseen parameters and/or geometries.
This is the \textit{online} stage, which consists of the following steps:
\begin{enumerate}
    \item Receive $\Omega$ and $\parVec$,
    \item Compute $\geoVec = \geoMapping(\Omega)$.
    \item For any $\pt_j$ for which the solution is needed, compute $\ptRef_j = \UCsys(\pt_j, \Omega)$.
    \item Evaluate $\sol(\pt_j; \parVec, \Omega) \simeq \NN(\ptRef_j, \parVec, \geoVec; \annParams)$. Typically, this operation can be vectorized to further increase the velocity of execution.
\end{enumerate}
We recall that both $\geoMapping$ and (if used) $\UCsys$ are defined case-by-case, depending on the application.

\section{Test cases} \label{sec:test-cases}

In this section, we present two test cases and provide details on the implementation choices we followed to apply the methods presented in Sec.~\ref{sec:methods}.

\subsection{Test Case 1: lid-driven cavity} \label{sec:test-cases:cavity}

Test Case 1 is based on the well-known stationary lid-driven cavity problem (see, e.g., \cite{botella1998benchmark}), for which we consider an extension with variable geometry.
The challenge of this test case is to capture the different vortex topologies formed for different Reynolds numbers and different aspect ratios of the geometry.


We consider a rectangular domain $\Omega_H = (0, 1) \times (0, H)$, with $H > 0$, and the following PDE (Navier-Stokes equations), where $\Gamma_H^D = \{ (x,y)^T \in \Omega_H \text{ s.t. } y = H \}$ denotes the top edge of the domain and $\Rey$ the Reynolds number:
\begin{equation} \label{eqn:cavity}
    \left\{
        \begin{aligned}
        &- \frac{1}{\Rey} \Delta \vel + ( \vel \cdot \nabla) \vel + \nabla \prs = \mathbf{0}
        && \quad \text{in }\Omega_H,  \\
        &\nabla \cdot  \vel = 0
        && \quad \text{in } \Omega_H,  \\
        & \vel = (1,0)^T
        && \quad \text{on } \Gamma_H^D, \\
        & \vel = \mathbf{0}
        && \quad \text{on } \partial\Omega_H \setminus \Gamma_H^D. \\
        \end{aligned}
    \right.
\end{equation}
The unknowns on the problem are the fluid velocity $\vel$ and pressure $\prs$.
The goal of Test Case 1 is to build an emulator that approximates the fluid velocity $\vel$, given the geometry $\Omega_H$ and the Reynolds number.
More precisely, we consider geometries with height $H$ in the interval $[0.5, 2]$:
\begin{equation*}
    \geoSpace = \{\Omega_H = (0, 1) \times (0, H) \text{, for } 0.5 \leq H \leq 2\}.
\end{equation*}
The physical parameter consists of $\parVec = \Rey$ and ranges in the interval $\parSpace = [10^2, 10^4]$.

\subsubsection{Training data generation}

To generate training data, we consider a Taylor-Hood Finite Element approximation of problem \eqref{eqn:cavity}.
We employ structured triangular computational meshes with a uniform space resolution of $h = 10^{-2}$.
We remark that, as a consequence, Finite Element approximations associated with different domains of $\geoSpace$ might feature different numbers of degrees of freedom.
To tackle Newton convergence issues for large $\Rey$, we equip the solver with an adaptive continuation ramp with respect to the Dirichlet datum.

To explore the set $\geoSpace \times \parSpace$, we employ a Monte Carlo approach by independently sampling from a uniform distribution for $H$ and a log-uniform distribution for $\Rey$.
After each FOM resolution, we export the velocity $\vel(\pt)$ at a set of points randomly selected within the domain $\Omega_H$.

\subsubsection{Geometrical landmarks}

Test Case 1 has an explicit parametrization of the domains in the set $\geoSpace$, the height $H$ being the parameter.
Henceforth, we define $\geoVec = H$ as the unique geometrical landmark, by setting $\geoMapping(\Omega_H) := H$.

\subsubsection{UC system}

A straightforward (and also effective) UC system for Test Case 1 consists in mapping each domain $\Omega_H \in \geoSpace$ into the unit square $\DomainRef := (0,1)^2$, through the transformation:
\begin{equation*}
    \hat{x} = x, \quad \hat{y} = y / H.
\end{equation*}
More precisely,
\begin{equation*}
    \UCsys\left( \begin{pmatrix} x \\ y\end{pmatrix}, \Omega_H \right) := \begin{pmatrix} x \\ y / H\end{pmatrix}.
\end{equation*}

\subsubsection{USM-Net architecture}

\begin{figure}[ht]
    \includegraphics[width=\columnwidth]{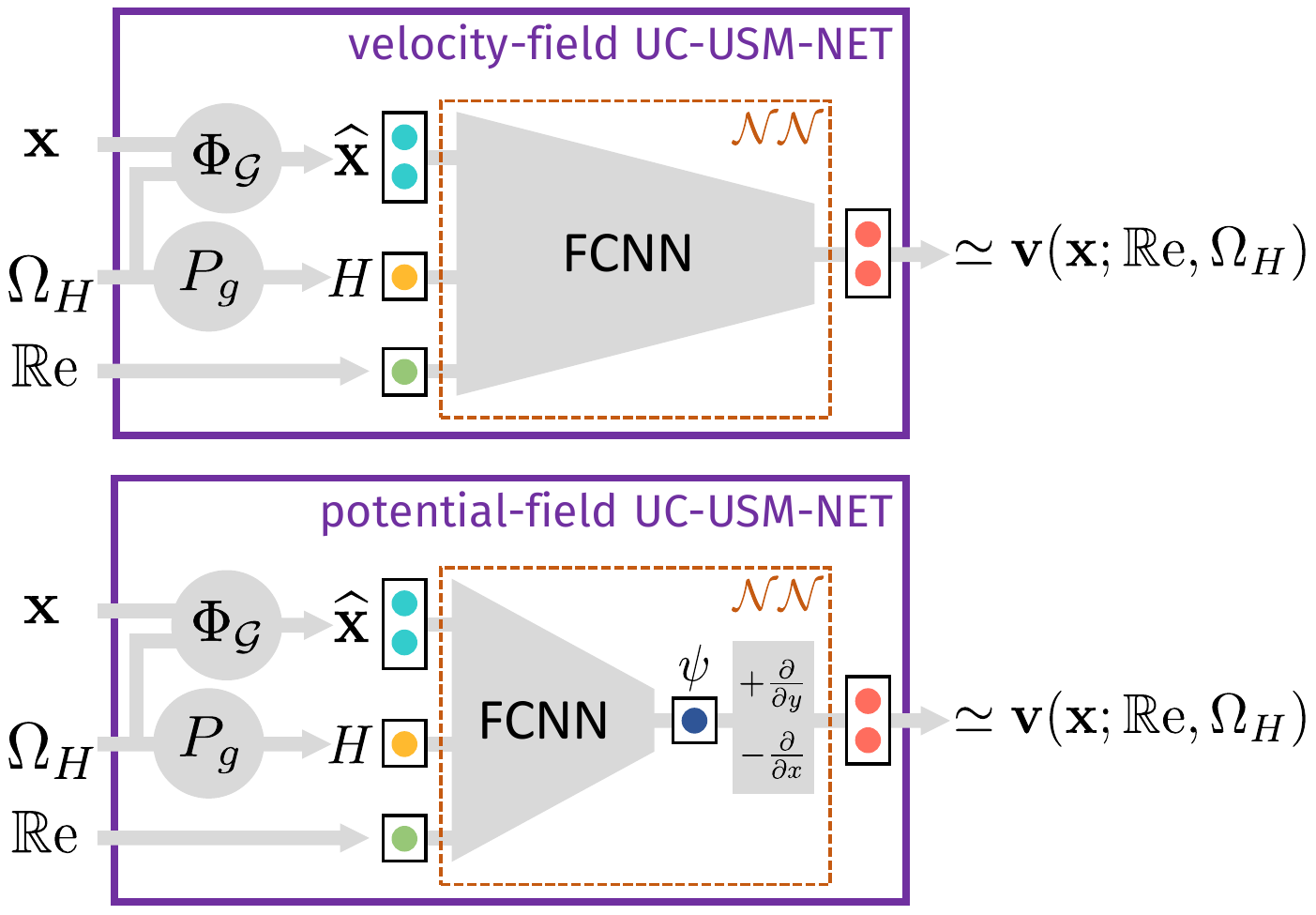}
    \caption{Test Case 1: comparison of velocity-field and potential-field architectures.}
    \label{fig:architecture_cavity}
\end{figure}

We consider two different ANN architectures to build USM-Nets for Test Case 1 (see Fig.~\ref{fig:architecture_cavity}).

\paragraph{Velocity-field architecture}
The first architecture for $\NN$ relies on a FCNN mapping $\pt$ (or $\ptRef$), $\parVec$ and $\geoVec$ into an approximation of $\vel(\pt; \Rey, \Omega_H)$.
To ease the FCNN training, we normalize both input and output data by mapping them in the interval $[-1, 1]$, and we preprocess the Reynolds number through a $\log$ transformation.
In conclusion, the FCNN features 4 input neurons and 2 output neurons.

\paragraph{Potential-field architecture}
As an alternative, we build a FCNN with a single output neuron, interpreted a the fluid flow potential $\pot(\pt; \Rey, \Omega_H)$, and we subsequently compute the approximation of the velocity field as:
\begin{equation} \label{eqn:potential}
    \vel(\pt; \Rey, \Omega_H) =
    \begin{pmatrix}
        +\frac{\partial}{\partial y} \pot(\pt; \Rey, \Omega_H) \\
        -\frac{\partial}{\partial x} \pot(\pt; \Rey, \Omega_H) \\
    \end{pmatrix}
\end{equation}
These operations are performed through Automatic Differentiation (AD) of the FCNN output.
We remark that we do not need a FOM-based potential $\pot$ for training data.
The training is done directly with respect to the velocity data.
The operations of \eqref{eqn:potential} represent indeed the last layer of the architecture $\NN$.
Input and outputs are normalized as for the velocity-field architecture.

\subsubsection{Loss function}

Since the goal of Test Case 1 is to reconstruct the velocity field with a focus on the vortex structure of the solution, we employ a discrepancy metric $\dist$ that emphasizes the role of flow direction at each point in the domain, including those with low flow intensity.
Specifically, we define
\begin{equation*}
    \begin{split}
        \dist(\sol, \tilde{\sol})
        =
        \| \sol -\tilde{\sol}\|^2
        +
        \left\| \frac{\sol}{\epsilon + \| \sol \| } - \frac{\tilde{\sol}}{\epsilon + \| \tilde{\sol} \|}  \right\|^2
    \end{split}
\end{equation*}
where $\epsilon = 10^{-4}$ is a small constant to avoid singularities.
The second term drives the USM-Net to accurately match the direction of the velocity.
Without this term, indeed, the flow direction would not be captured well in the regions with low flow magnitude, due to the low contribution in the first term of the loss.

Moreover, we augment the loss function with the following physics-based regularization term, aimed at enforcing the satisfaction of the Dirichlet boundary conditions:
\begin{equation*}
\begin{array}{l}
    \RegTerm(\annParams) = \frac{1}{N_{\text{BC}} \numPoints{\text{BC}}} \sum_{i=1}^{N_{\text{BC}}} \sum_{j=1}^{\numPoints{\text{BC}}}
    \| \mathbf{u}_{\text{BC},j}^i  - \vel_{\text{BC}}(\pt_{\text{BC},j}^i) \|^2, \\
    \mathbf{u}_{\text{BC},j}^i = \NN(\pt_{\text{BC},j}^i, \parVec^{\text{BC}, i}, \geoVec^{\text{BC}, i}; \annParams)
\end{array}
\end{equation*}
where $\parVec^{\text{BC}, i} \in [10^{2}, 10^{4}]$ and $\geoVec^{\text{BC}, i} \in [0.5, 2]$ are sample points, $\pt_{\text{BC},j}^i$ is a set of points belonging to the boundary, and $\vel_{\text{BC}}$ is the Dirichlet datum (see \eqref{eqn:cavity}).

\subsection{Test Case 2: coronary bifurcation} \label{sec:test-cases:coronary}

As a second test case, we consider the problem of predicting blood flow and pressure field within a coronary bifurcation in the presence of stenosis.

\begin{figure}[ht]
    \centering
    \includegraphics[width = \columnwidth]{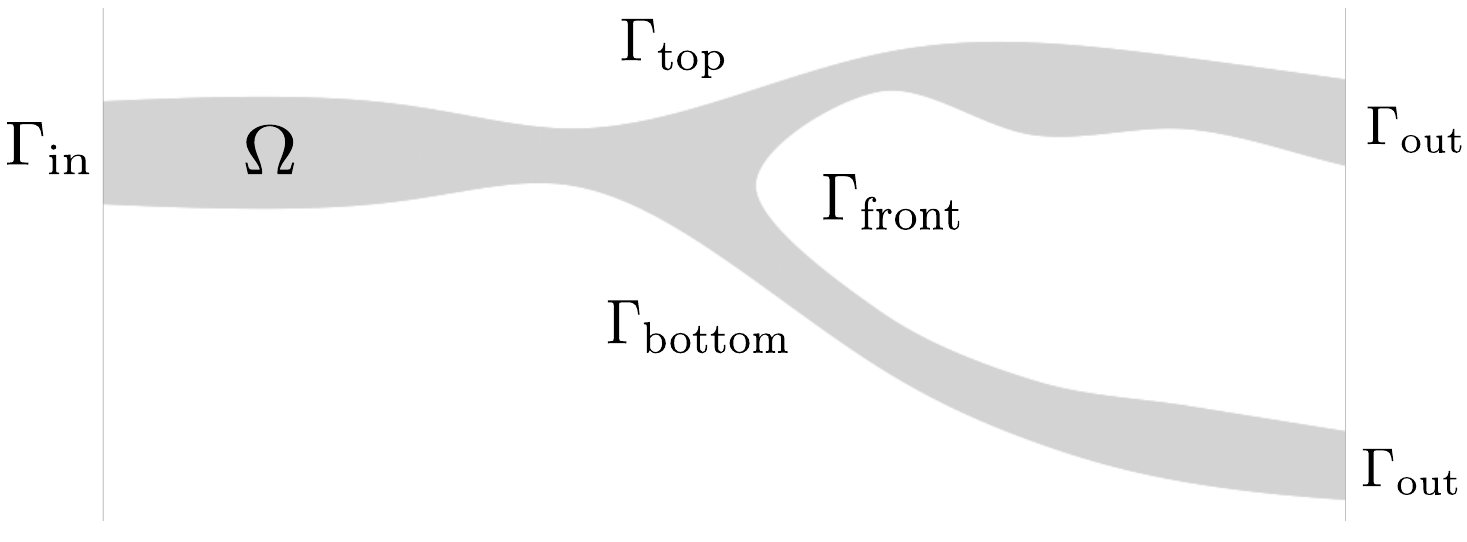}\\
    \caption{Test Case 2: example of computational domain and corresponding boundaries.}
    \label{fig:coro_domain}
\end{figure}

More precisely, we consider a computational domain $\Omega$, corresponding to the 2D representation of a section of a coronary artery with a bifurcation.
We denote by $\Gamma_{\text{in}}$ the portion of the boundary corresponding to the inlet, by $\Gamma_{\text{out}}$ the two outlets and by $\Gamma_{\text{wall}} = \Gamma_{\text{top}} \cup \Gamma_{\text{bottom}} \cup \Gamma_{\text{front}}$ the vessel wall.
In this test case, we will consider many different computational domains, each representing a coronary bifurcation in a different virtual patient.
An example domain is represented in Fig.~\ref{fig:coro_domain}.

\begin{figure*}
    \centering
    \includegraphics[]{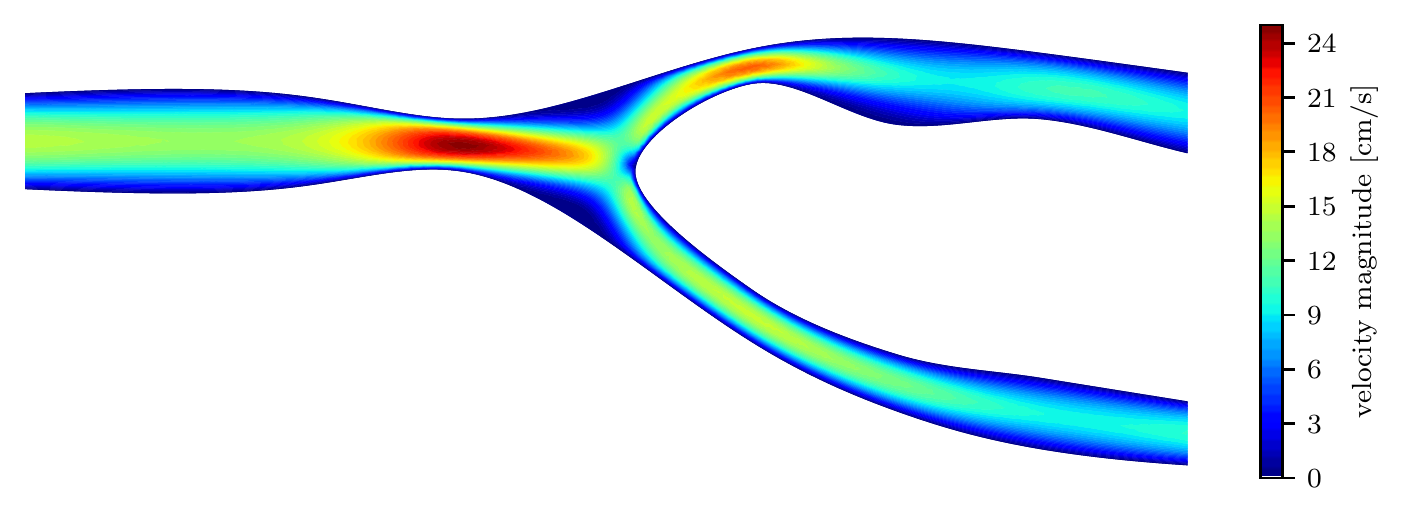}\\
    \includegraphics[]{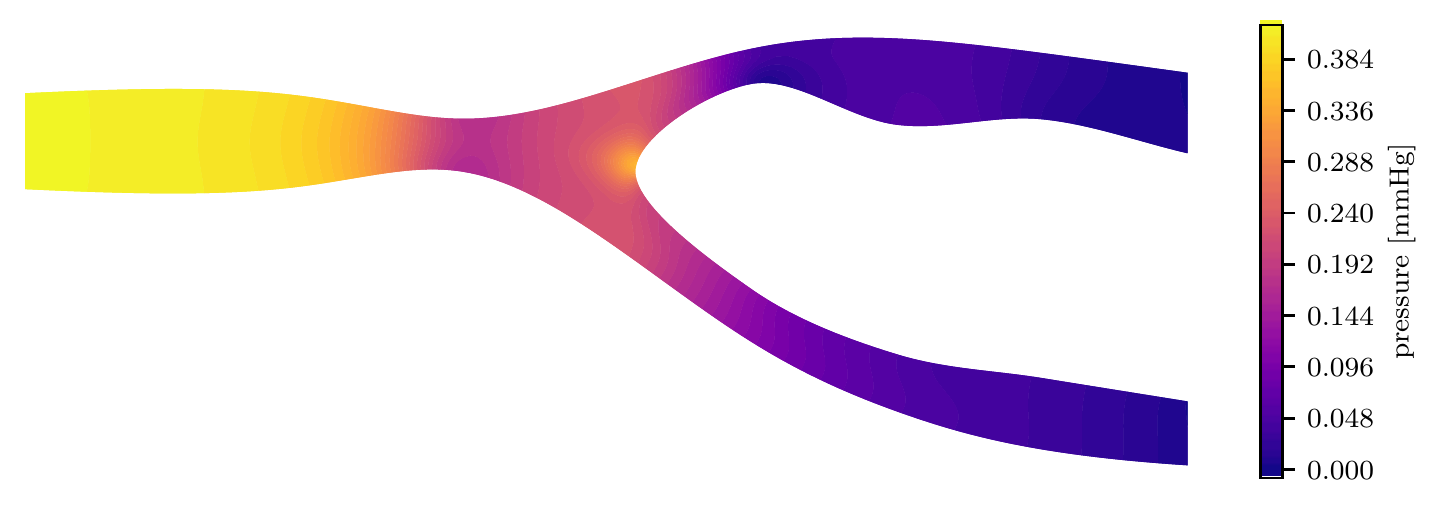}\\
    \caption{Test Case 2: numerical solution of \eqref{eqn:coronary} on the computational domain of Fig.~\ref{fig:coro_domain}. Top: velocity field; bottom: pressure field.}
    \label{fig:coro_example_solution}
\end{figure*}

We consider the following stationary Navier-Stokes model, describing the steady-state fluid flow in the coronary bifurcation:
\begin{equation} \label{eqn:coronary}
    \left\{
        \begin{aligned}
        &- \nu \Delta \vel +  ( \vel \cdot \nabla) \vel + \frac{1}{\rho} \nabla \prs = \mathbf{0}
        && \quad \text{in }\Omega,  \\
        &\operatorname{div}  \vel = 0
        && \quad \text{in } \Omega,  \\
        & \vel = \vel_{\text{in}}
        && \quad \text{on } \Gamma_{\text{in}}, \\
        & \vel = \mathbf{0}
        && \quad \text{on } \Gamma_{\text{wall}}, \\
        & \nu \frac{\partial \vel}{\partial \mathbf{n}} - \frac{1}{\rho} \prs \, \mathbf{n} = \mathbf{0}
        && \quad \text{on } \Gamma_{\text{out}}, \\
        \end{aligned}
    \right.
\end{equation}
where $\nu = \SI{4.72}{\milli\liter\squared\per\second}$ is the kinematic viscosity of blood and $\rho = \SI{1060}{\kilogram\per\meter\cubed}$ its density.
At the inlet, we consider a parabolic profile with a maximum velocity equal to $\SI{14}{\centi\meter\per\second}$.
In Fig.~\ref{fig:coro_example_solution}, we show the numerical solution of \eqref{eqn:coronary} in the example computational domain of Fig.~\ref{fig:coro_domain}.

\subsubsection{Geometrical variability and landmarks}

\begin{figure*}
    \centering
    \includegraphics[width=\textwidth]{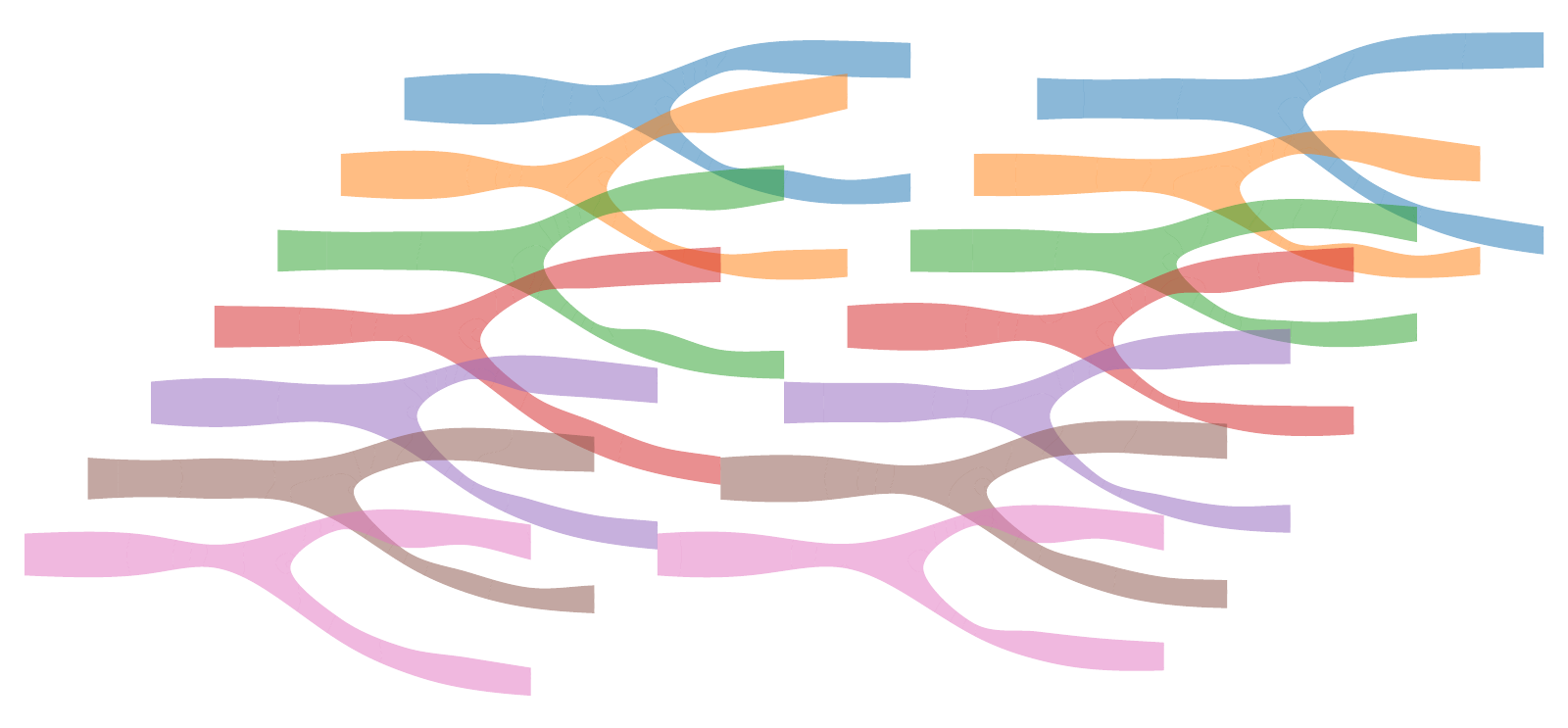}\\
    \caption{Test Case 2: representation of some of the geometries $\Omega \in \geoSpace$ included in the training dataset.}
    \label{fig:coro_geometries}
\end{figure*}

The aim of Test Case 2 is to build a reduced model predicting the pressure and velocity fields in an arbitrary domain representing a coronary bifurcation.
We synthetically generate a large number of different computational domains corresponding to many virtual patients.
To do this, we use splines obtained by randomly varying their parameters in suitable intervals, defined to reflect the realistic variability observed among patients \cite{chiastra2016coronary}.
Notice that the geometries thus obtained may present stenoses, of a more or less acute degree, located upstream of the bifurcation or in the two branches downstream of it.
A subsample of the geometries obtained following this procedure is displayed in Fig.~\ref{fig:coro_geometries}.

\begin{figure}[ht]
    \centering
    \includegraphics[width = \columnwidth]{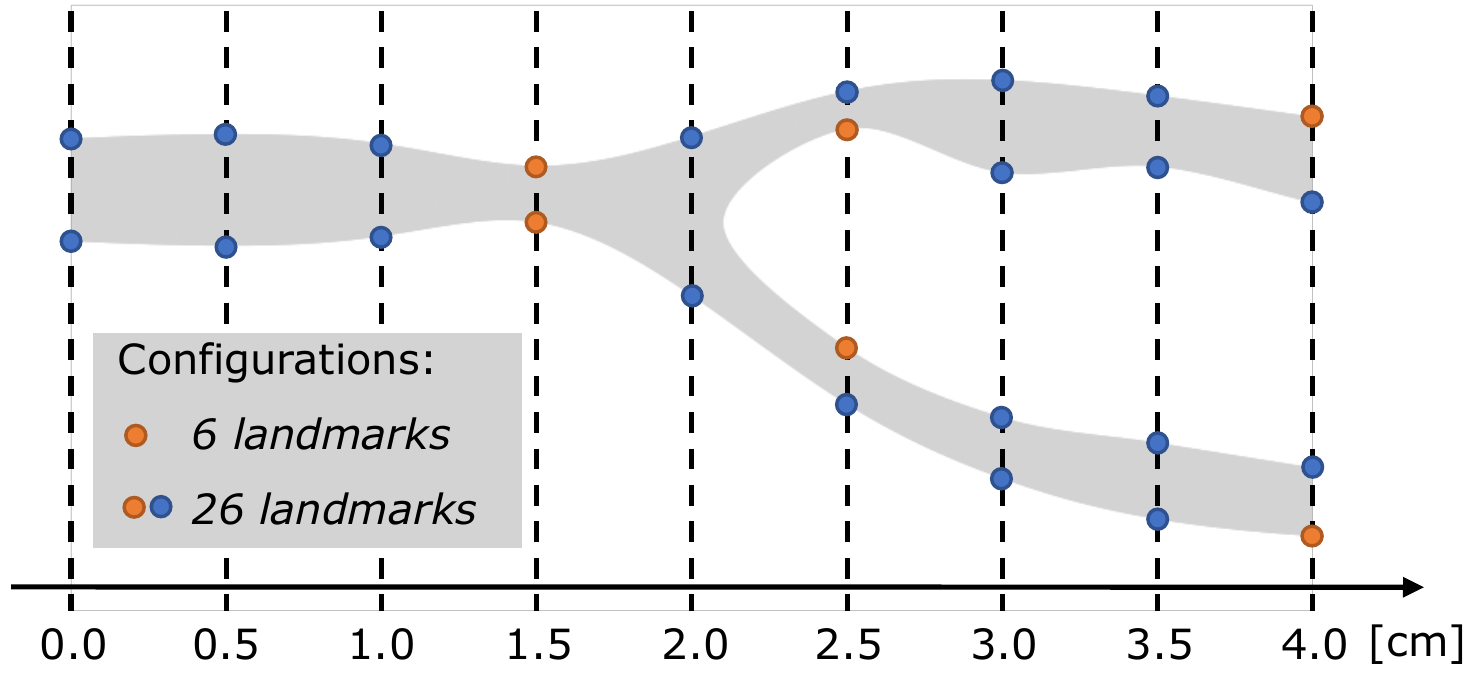}\\
    \caption{Test Case 2: geometrical landmarks $\geoVec$.}
    \label{fig:coro_landmarks}
\end{figure}

Due to the lack of an explicit parameterization of these geometries (a common issue when dealing with domains from real patients), we need to define geometrical landmarks to characterize each geometry.
To this end, we use an operatively light procedure that can also be easily adopted in a clinical context.
Specifically, we define landmarks as the coordinates $y$ of the vessel wall corresponding to a set of predefined coordinates $x$.
These coordinates contain information regarding the lumen thickness at various levels and the possible presence of stenosis.
Note that, in clinical practice, these landmarks can be easily derived directly from medical imaging, without the need to construct a computational mesh.
In this test case, we will consider two sets of landmarks containing, respectively, 26 and 6 coordinates (see Fig.~\ref{fig:coro_landmarks}).
Clearly, the first set is much more informative than the second one.
The aim is to test the robustness of the proposed methods in the case where the landmarks provide a poor description of the geometry, and are not able to exhaustively capture its variability.

\subsubsection{UC system}

\begin{figure*}
    \centering
    \includegraphics[]{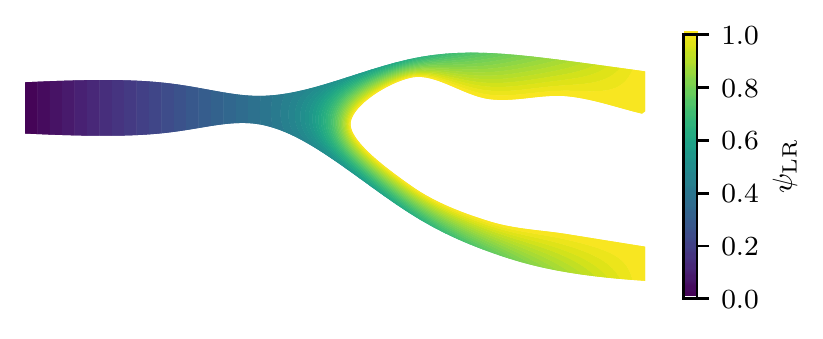}
    \includegraphics[]{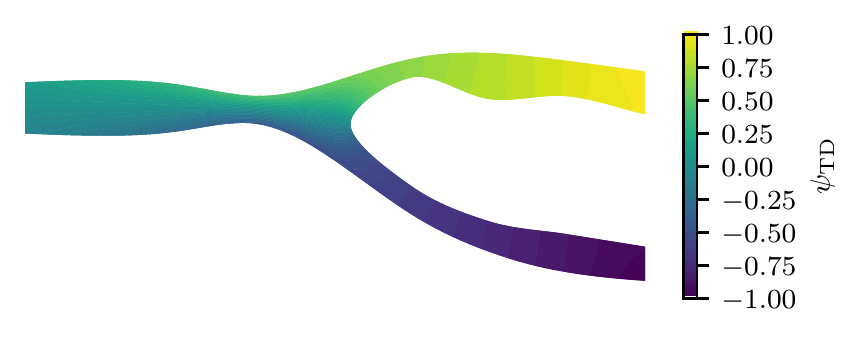}
    \caption{Test Case 2: fields $\UCcoroLR$ and $\UCcoroTD$ associated with the domain of Fig.~\ref{fig:coro_example_psi}.}
    \label{fig:coro_example_psi}
\end{figure*}

Differently than in Test Case 1, where, thanks to the simplicity of the domains of the $\geoSpace$ space, it was possible to explicitly define a UC system, in Test Case 2 the construction of a UC system is not an immediate task.
We propose to rely on two Laplacian-based fields, which define the inlet-outlet and top-bottom directions, respectively.
More precisely, given a geometry $\Omega \in \geoSpace$, we define the fields $\UCcoroLR\colon\Omega \to \DomainRef$ and $\UCcoroTD\colon\Omega \to \DomainRef$ as the solutions of the following differential problems:
\begin{equation} \label{eqn:UCcoroLR}
    \left\{
        \begin{aligned}
        &- \Delta \UCcoroLR = 0
        && \quad \text{in }\Omega,  \\
        & \UCcoroLR = 0
        && \quad \text{on } \Gamma_{\text{in}}, \\
        & \UCcoroLR = 1
        && \quad \text{on } \Gamma_{\text{out}} \cup \Gamma_{\text{front}}, \\
        & \UCcoroLR = \frac{x - x_{\min}}{x_{\max} - x_{\min}}
        && \quad \text{on } \Gamma_{\text{top}} \cup \Gamma_{\text{bottom}}, \\
        \end{aligned}
    \right.
\end{equation}
\begin{equation} \label{eqn:UCcoroTD}
    \left\{
        \begin{aligned}
        &- \Delta \UCcoroTD = 0
        && \; \text{in }\Omega,  \\
        & \UCcoroTD = + \alpha + (1 - \alpha)\frac{x - x_{\min}}{x_{\max} - x_{\min}}
        && \; \text{on } \Gamma_{\text{top}}, \\
        & \UCcoroTD = -\alpha - (1 - \alpha)\frac{x - x_{\min}}{x_{\max} - x_{\min}}
        && \; \text{on } \Gamma_{\text{bottom}}, \\
        & \frac{\partial \UCcoroTD}{\partial \mathbf{n}} = 0
        && \; \text{on } \Gamma_{\text{in}} \cup \Gamma_{\text{out}} \cup \Gamma_{\text{front}}. \\
        \end{aligned}
    \right.
\end{equation}
In Fig.~\ref{fig:coro_example_psi} we show the fields $\UCcoroLR$ and $\UCcoroTD$ obtained for the domain of Fig.~\ref{fig:coro_domain}.
The field $\UCcoroLR$ bridges the inlet region (i.e. $\Gamma_{\text{in}}$, where $\UCcoroLR = 0$) with the frontal region of the domain (i.e. $\Gamma_{\text{out}} \cup \Gamma_{\text{front}}$, where $\UCcoroLR = 1$).
Conversely, $\UCcoroTD$ defines the proximity of each lumen point to the upper wall relative to the lower wall.
Setting a constant $\alpha < 1$ allows better differentiation of the $\UCcoroTD$ field within each branch downstream of the bifurcation.
Specifically, we set $\alpha = 0.1$.

\begin{figure*}
    \centering
    \includegraphics[width = 0.8\textwidth]{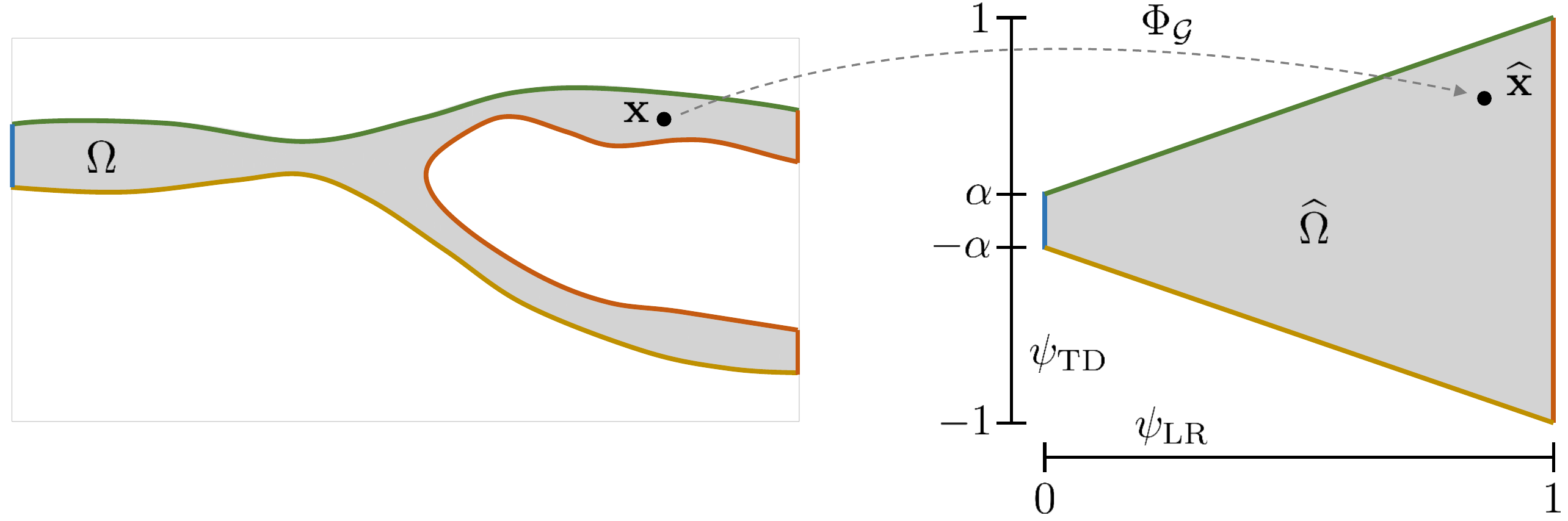}
    \caption{Test Case 2: representation of the UC system $\UCsys$. Left: physical domain $\Omega$; right: reference domain $\DomainRef$.}
    \label{fig:coro_UC}
\end{figure*}

The UC system is thus defined as:
\begin{equation*}
    \ptRef =
    \UCsys( \pt, \Omega ) :=
    \begin{pmatrix}
        \UCcoroLR(\pt; \Omega) \\
        \UCcoroTD(\pt; \Omega)
    \end{pmatrix}.
\end{equation*}
In Fig.~\ref{fig:coro_UC} we show the reference domain $\DomainRef$ and the mutual correspondences between the boundary of the physical and reference domains.

\subsubsection{Training data generation}

To generate training data, we employ the Finite Element solver described for Test Case 1 (see Sec.~\ref{sec:test-cases:cavity}).
For space discretization, we consider triangular computational meshes with a space resolution of nearly $h = \SI{0.2}{\milli\meter}$.
The UC coordinates are obtained by solving the differential problems \eqref{eqn:UCcoroLR} and \eqref{eqn:UCcoroTD} with $P1$ Finite Elements on the same computational mesh.

\subsubsection{USM-Net architecture}

In Test Case 2, we employ a FCNN connecting $\pt$ or $\ptRef$ (depending on whether a PC-USM-Net or a UC-USM-Net is used), $\parVec$ and $\geoVec$ into an approximation of $\sol(\pt; \Omega_H)$ where the solution $\sol = (\vel, \prs)$ is given by the pair velocity-pressure.
Similarly to Test Case 2, we add a normalization layer before the input and after the output layers of the FCNN, to constrain each input and each output in the interval $[-1, 1]$.

\subsubsection{Loss function}
We employ a purely black-box loss function, as defined in \eqref{eqn:loss}, with a quadratic discrepancy metric
$$
\dist(\sol, \tilde{\sol}) = \| \sol -\tilde{\sol}\|^2
$$
and without any regularization terms.

\section{Results} \label{sec:results}

In this section, we present the results obtained by applying the methods presented in Sec.~\ref{sec:methods} to the two test cases of Sec.~\ref{sec:test-cases}.
These results have been obtained using TensorFlow \cite{tensorflow2015whitepaper} and the optimizers of SciPy \cite{SciPy2020}.

\subsection{Test Case 1: lid-driven cavity}

\begin{figure*}
    \centering
    \includegraphics[width=0.9\textwidth]{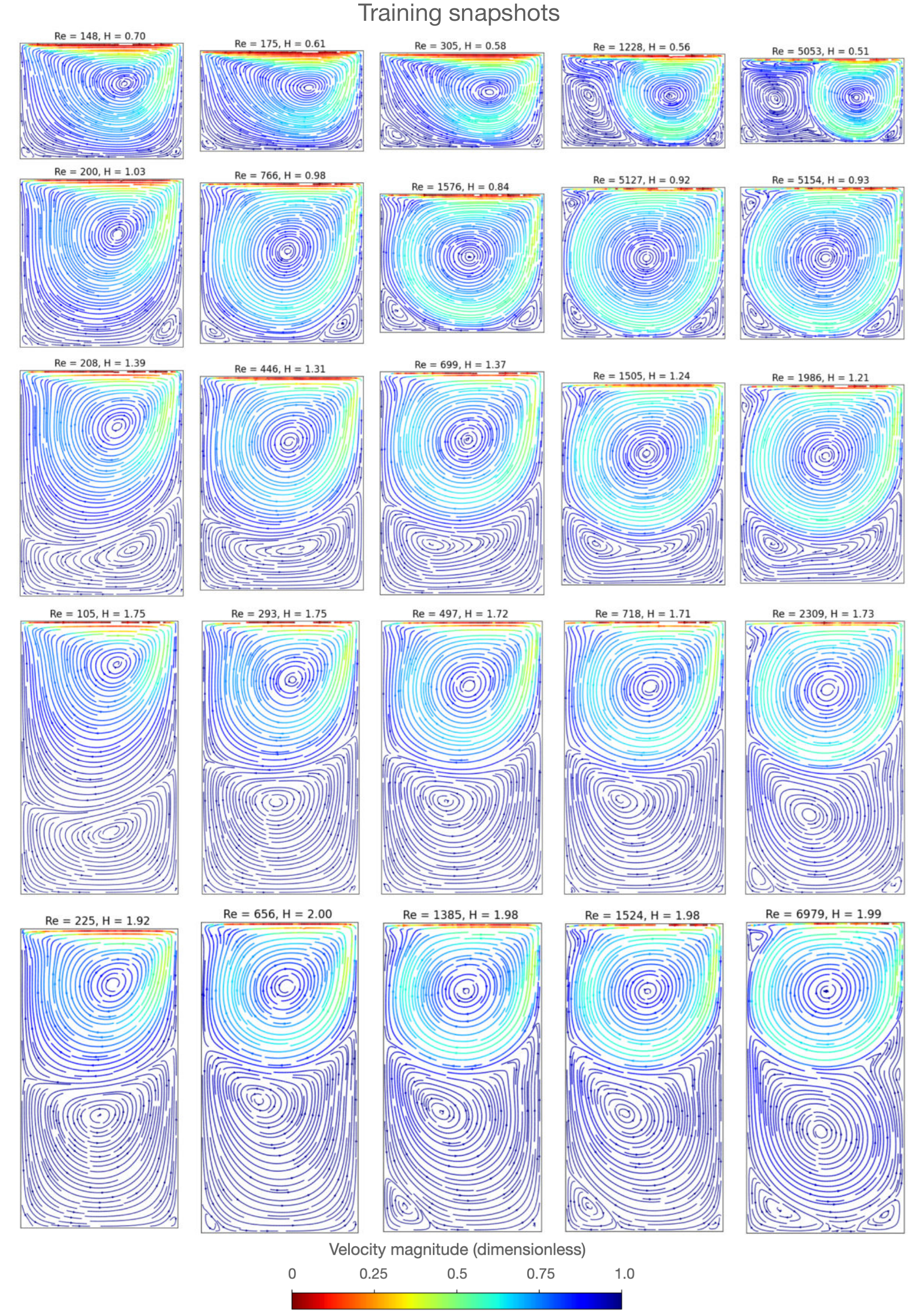}\\
    \caption{
    Test Case 1: comparison of some numerical solutions that constitute the training set. The dataset is generated by approximating the FOM \eqref{eqn:cavity} for random samples of the values of the physical and geometrical parameters.}
    \label{fig:cavity_comparison_train}
\end{figure*}

We construct a training set of $\numTrain = 400$ numerical simulations obtained by randomly sampling the height $H$ and the physical parameter $\parVec = \Rey$. Some examples of streamlines resulting from numerical solutions by varying the parameters are reported in Fig.~\ref{fig:cavity_comparison_train}. For each geometry, we subsample the solution in $\numPoints{} = 360$ random points.

\begin{figure*}
    \centering
    \includegraphics[width=1.0\textwidth]{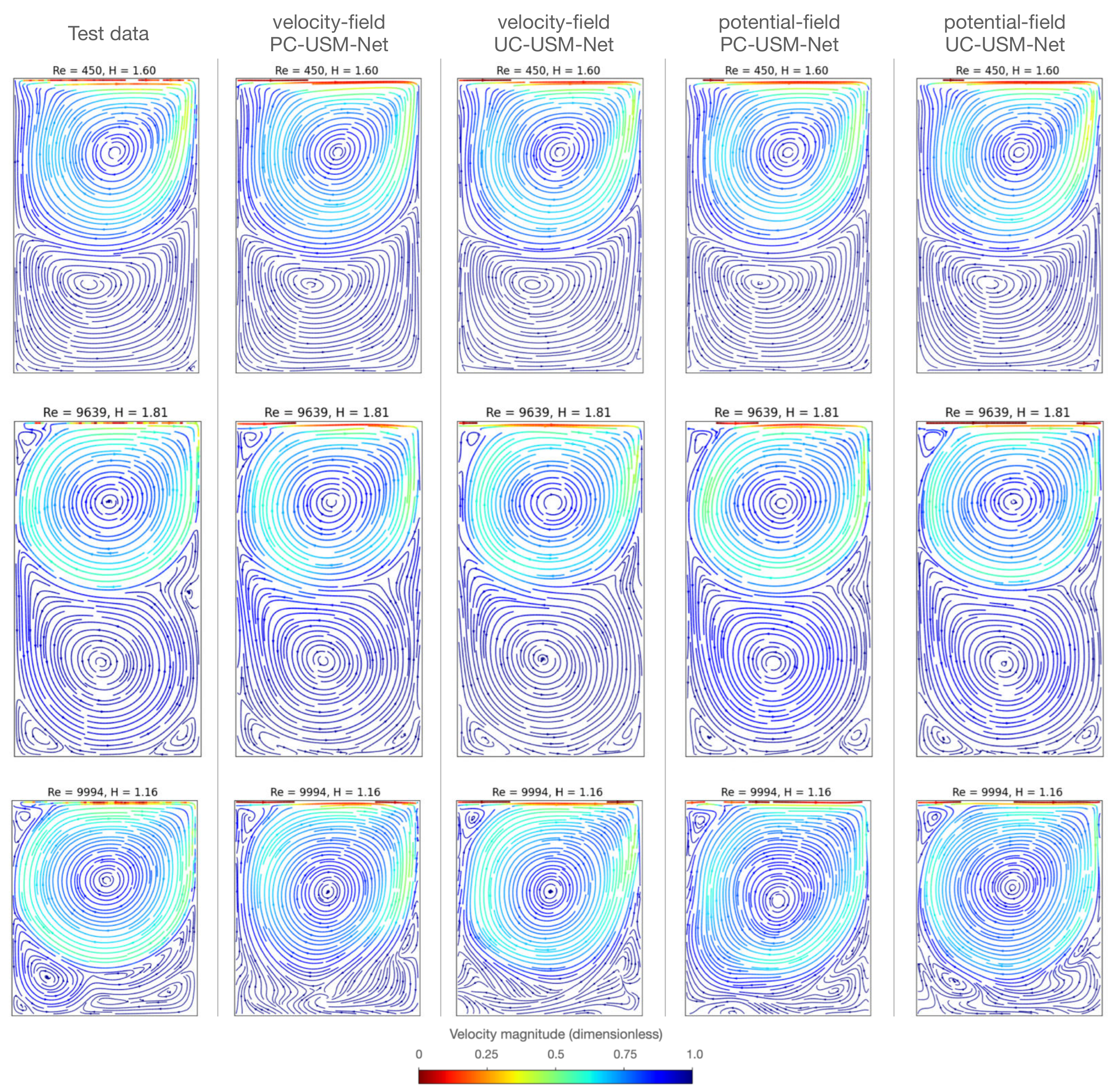}\\
    \caption{
    Test Case 1: comparison of ANN streamlines reconstruction on three numerical solutions from the test set. Test cases are selected to display the range of ANN reconstruction errors in the test set: from minor (best case scenario: first row) to significant (worst case scenario: last row).}
    \label{fig:Cavity_strategy_comparison}
\end{figure*}

By varying the dimensions of this dataset, we train several USM-Nets, formed by FCNNs made of 3 inner layers, respectively consisting of 30, 20, and 10 neurons.
We consider four configurations:
\begin{enumerate}
\item velocity-field {PC-USM-Net}: receiving as inputs the two spatial coordinates, the physical and the geometrical parameters, and producing as outputs the two velocity components;
\item velocity-field {UC-USM-Net}: receiving as input the two universal coordinates, the physical and the geometrical parameters, and producing as outputs the two velocity components;
\item potential-field {PC-USM-Net}: receiving as inputs the two spatial coordinates, the physical and the geometrical parameters, and producing as output the two velocity components computed from the fluid flow potential;
\item potential-field {UC-USM-Net}: receiving as inputs the two universal coordinates, the physical and the geometrical parameters, and producing as output the two velocity components computed from the fluid flow potential.
\end{enumerate}
For each configuration, we perform $500$ epochs of the ADAM optimizers \cite{kingma2014adam} followed by $20000$ epochs of the BFGS method \cite{goodfellow2016deep} to ensure convergence of the optimizer.
For the case of a training set composed of $\numTrain = 100$ numerical simulations, we post-process the velocity field to display streamlines. In Fig.~\ref{fig:Cavity_strategy_comparison} we report the streamlines resulting from the different ANN configurations for three test cases extracted from the $40$ numerical simulations that formed the test set: they represent the best, the average, and the worst-case scenarios, respectively.

\begin{figure*}
    \centering
    \includegraphics[width=1.0\textwidth]{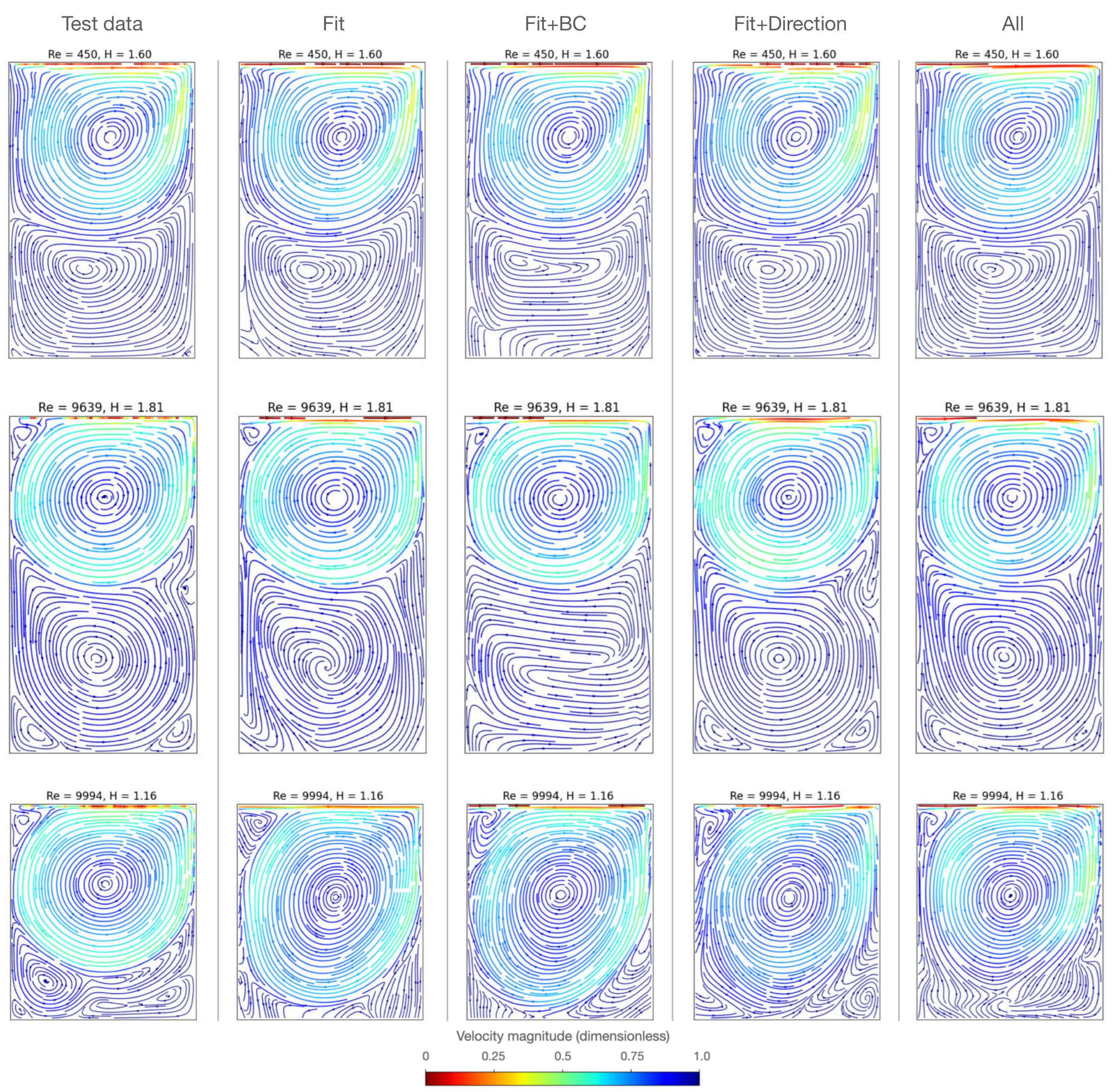}\\
    \caption{
    Test Case 1: comparison of velocity-field PC-USM-Net streamlines reconstruction on three numerical solutions from the test set on varying the definition of the loss function. We consider a loss function composed of the data misfit term (second column), combined with boundary conditions (third column) or with the direction regularization (fourth column). Test cases are selected to display the range of ANN reconstruction errors in the test set: from minor (best case scenario: first row) to significant (worst case scenario: last row).}
    \label{fig:Cavity_loss_comparison_model}
\end{figure*}

\begin{figure*}
    \centering
    \includegraphics[width=1.0\textwidth]{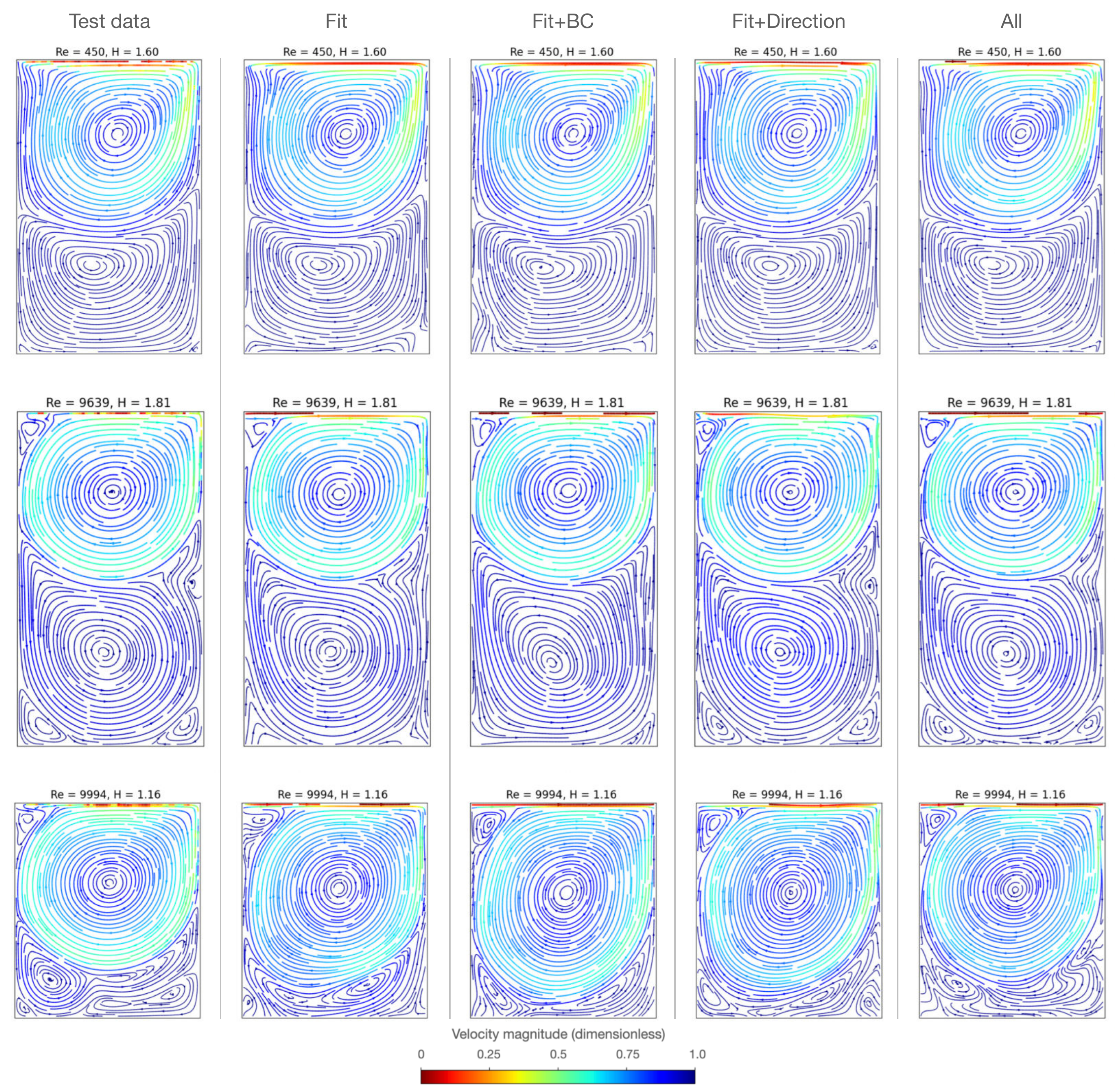}\\
    \caption{
    Test Case 1: comparison of potential-field UC-USM-Net streamlines reconstruction on three numerical solutions from the test set on varying the definition of the loss function. We consider a loss function composed of the data misfit term (second column), combined with boundary conditions (third column) or with the direction regularization (fourth column). Test cases are selected in order to display the range of ANN reconstruction errors in the test set: from minor (best case scenario: first row) to significant (worst case scenario: last row).}
    \label{fig:Cavity_loss_comparison_potential}
\end{figure*}

Training the ANN with a loss function composed only of the data misfit term results insufficient for an accurate reconstruction of the streamlines, especially in low-velocity areas. In Figs.~\ref{fig:Cavity_loss_comparison_model} and \ref{fig:Cavity_loss_comparison_potential}, we show the effect of the loss function components described in Section \ref{sec:test-cases:cavity} on the streamlines reconstruction.

\begin{figure*}
    \centering
    \includegraphics[width=1.0\textwidth]{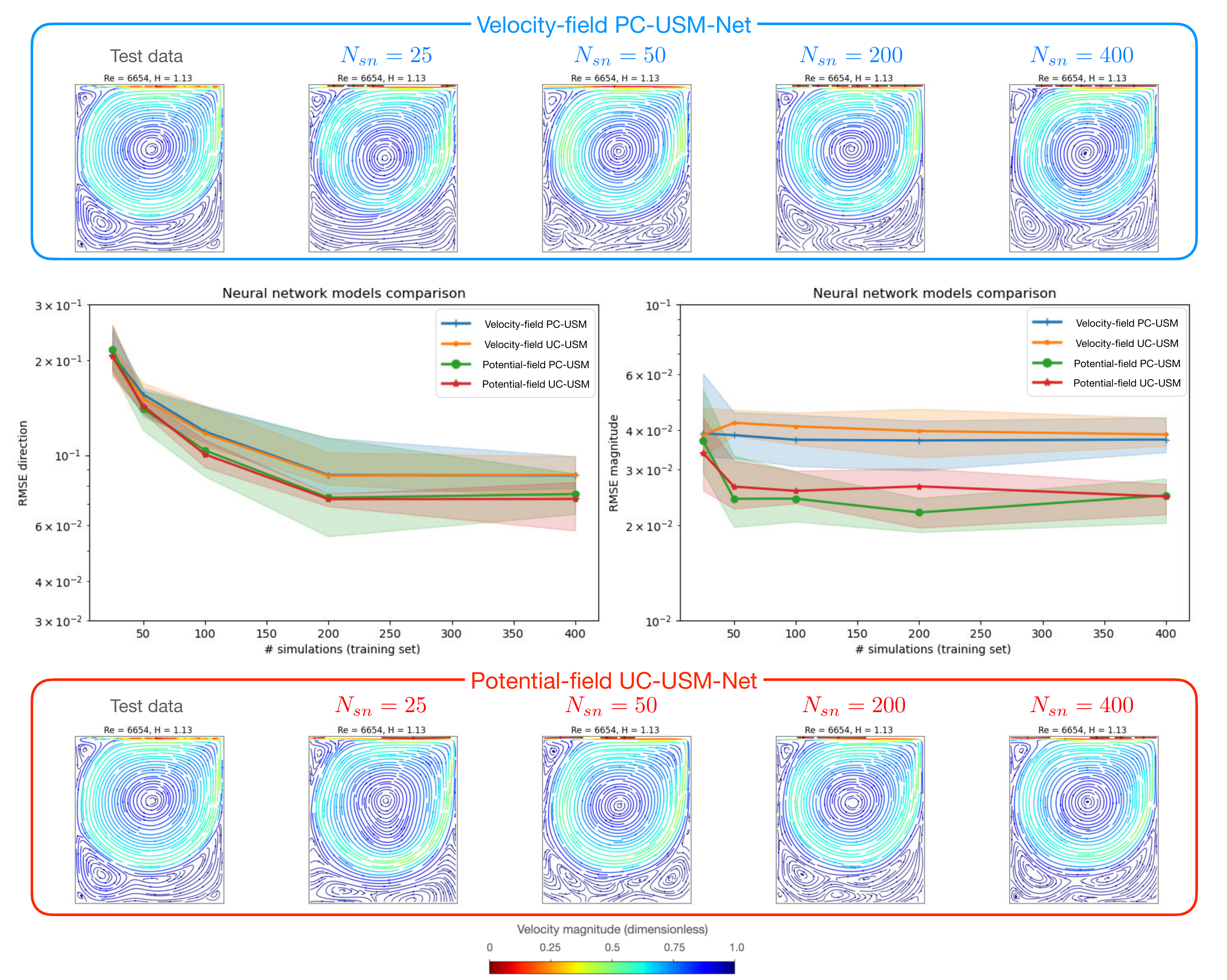}\\
    \caption{
    Test Case 1: RMSE on the velocity magnitude and direction on the test set, made of $40$ numerical simulations sampled in $10000$ points. We compare the four different configurations of ANN on varying the dimension of the training set. The error in the direction is significantly larger than the error in the magnitude and is inversely proportional to the number of simulations of the training set. }
    \label{fig:Cavity_RMSE}
\end{figure*}

To assess the generalization error on the test set resulting from different training configurations, we repeat the training considering only a subset of $[25,50,100,200,400]$ numerical simulations as training set.
The test set comprises $40$ numerical simulations sampled in $10000$ points.
We repeat the training with 10 different random initializations of the ANN weights and biases.
The average root mean squared error (RMSE) of the velocity magnitude and direction are reported in Fig.~\ref{fig:Cavity_RMSE}, together with bands indicating the maximum and the minimum values for the different training set dimensions and ANN configurations.

\subsection{Test Case 2: coronary bifurcation}

We consider a training set consisting of $\numTrain = 500$ different geometries.
For each geometry, we take the solution in $\numPoints{} = 1000$ randomly generated points.
Then, for both the landmark configurations (26 landmarks and 6 landmarks), we train a PC-USM-Net and a UC-USM-Net.
We consider a FCNN with 4 inner layers, respectively consisting of 20, 15, 10, and 5 neurons.
This architecture has been tuned in order to minimize the validation error on a set of 100 geometries not included in the training dataset.
To train the FCNN weights and biases, we run 200 iterations of the Adam optimizer \cite{kingma2014adam} and, subsequently, 5000 iterations of the BFGS algorithm \cite{goodfellow2016deep}.
We perform 10 different training runs for each configuration, starting from different random initializations of the FCNN parameters.
Each training run lasts about 45 minutes on a laptop equipped with Intel Core i7-1165G7 CPU (2.80 GHz).

\begin{figure*}
    \centering
    \includegraphics[]{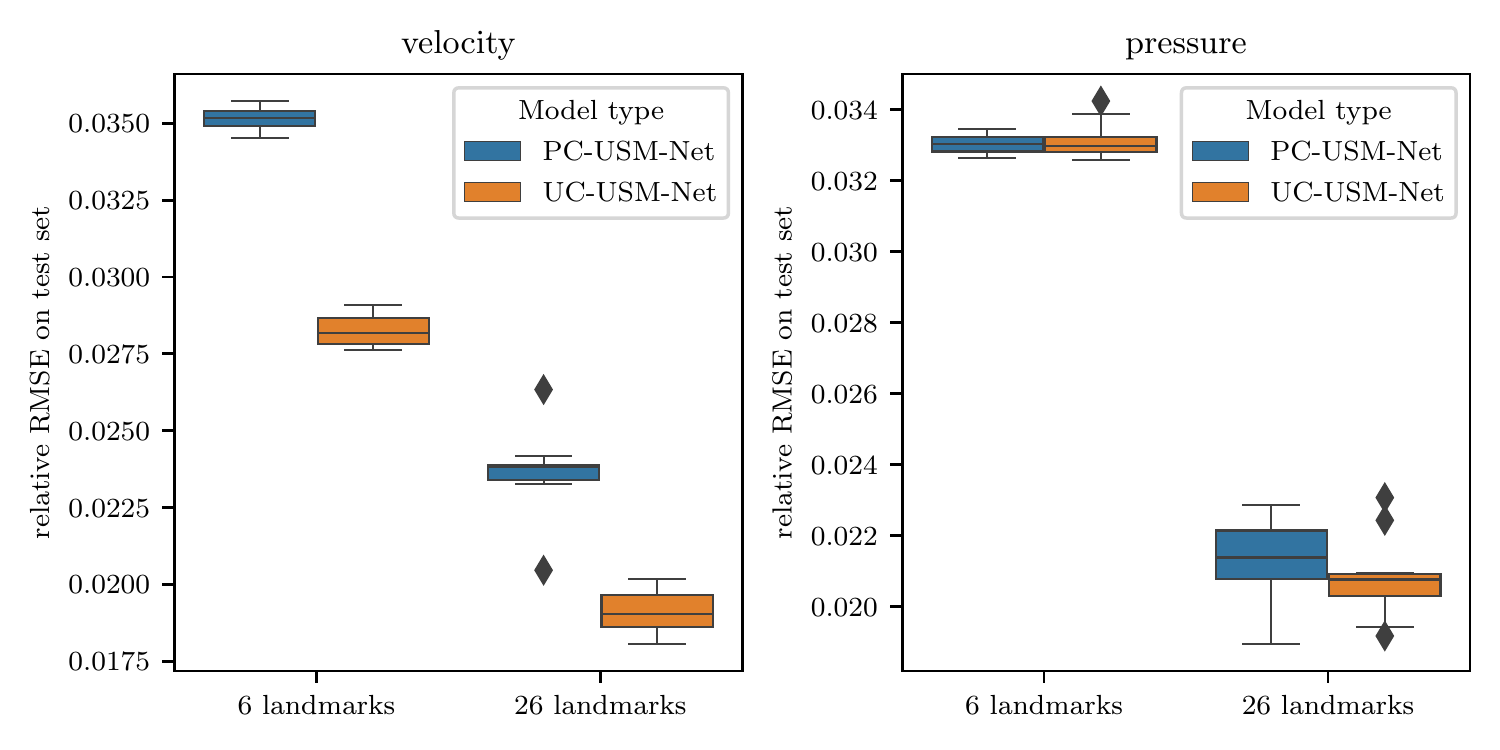}\\
    \caption{
    Test Case 2: boxplots of the errors on the test dataset obtained with 26 landmarks and 6 landmarks and with the PC-USM-Net and the UC-USM-Net architecture.
    The boxplots refer to 10 training runs obtained starting from different random initializations of the ANN weights and biases.
    Left: error on the velocity field; right: error on the pressure field.
    }
    \label{fig:coro_boxplot_up}
\end{figure*}

In Fig.~\ref{fig:coro_boxplot_up}, we show boxplots of the errors associated with a testing dataset of 100 geometries,  neither in the training nor in the validation dataset.
As expected, USM-Nets that are provided with 26 landmarks generate more accurate predictions than those that are aware of only 6 landmarks.
However, we notice that USM-Nets based on only 6 landmarks still have a noticeable accuracy (relative RMSE error of about 3\% on both the velocity and pressure).
This figure is compatible with the levels of precision typically required in clinical practice.

Furthermore, the boxplots show that the use of UC can significantly enhance the performance of the USM-Net.
The improvement is all the more evident in the case of the velocity field, compared to the pressure field.

\begin{figure*}
    \centering
    \includegraphics[width=\textwidth]{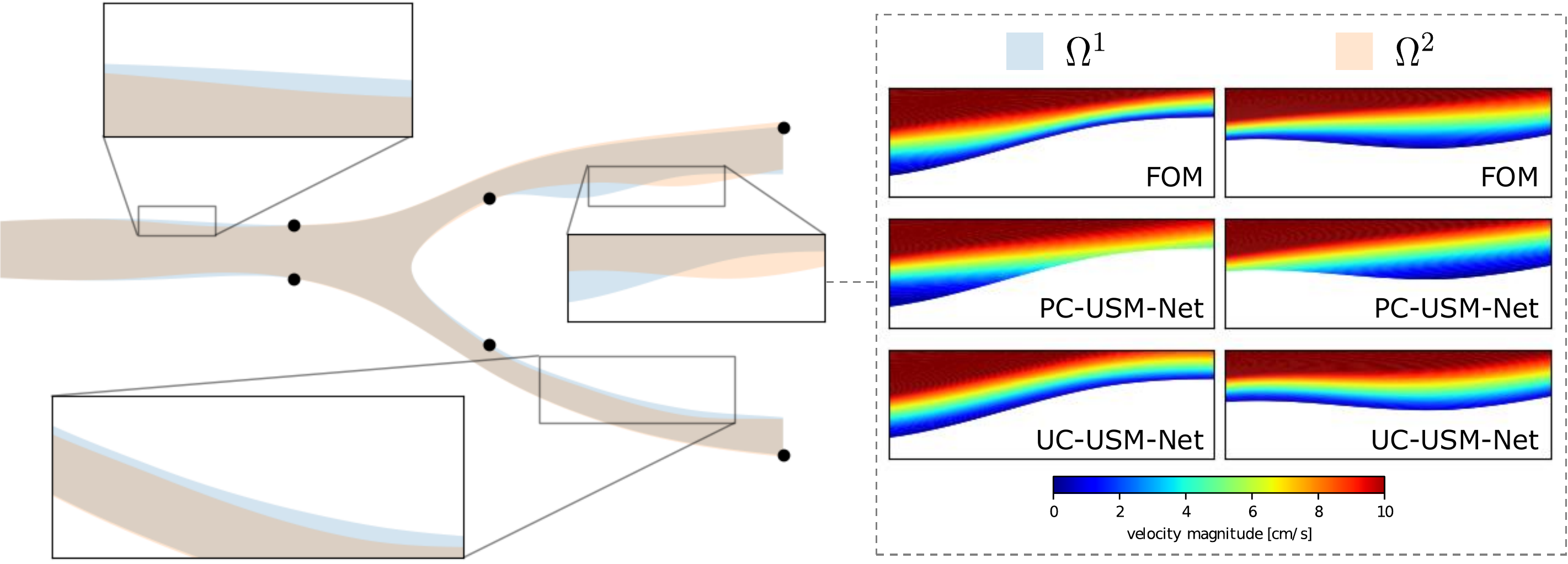}\\
    \caption{Test case 2: two domains ($\Omega^1$ and $\Omega^2$) belonging to the testing set that feature almost identical landmarks, that is $\geoMapping(\Omega^1) \simeq \geoMapping(\Omega^2)$ (left).
    On the right, a detail of the velocity field is compared among FOM solution, PC-USM-Net, and UC-USM-Net surrogates.
    The color map is intentionally flattened towards low values to highlight velocity variations near the edge.}
    \label{fig:coro_geo_variability}
\end{figure*}

In order to highlight the role that a UC system has in improving the generalization accuracy of USM-Nets, we consider the pair of domains in the testing set that are characterized by the most similar landmarks.
More precisely, these domains, which we will call $\Omega^1$ and $\Omega^2$, are such that $\|\geoMapping(\Omega^1) - \geoMapping(\Omega^2) \| < 10^{-4}$.
Since landmarks characterize the domain at some control points, two domains with very similar landmarks may differ significantly away from the control points.
This is what happens for the two domains considered, in particular in correspondence of the upper outflow track (see Fig.~\ref{fig:coro_geo_variability}, left).
On the right side of Fig.~\ref{fig:coro_geo_variability} we show a detail of the velocity field obtained for these two domains with a PC-USM-Net and a UC-USM-Net, in comparison with the reference solution obtained by means of the FOM.
We recall that the domain $\Omega$ affects the PC-USM-Net result only through the landmarks $\geoVec$.
Therefore, the PC-USM-Net will provide the same solution for two geometries with identical landmarks.
As shown by Fig.~\ref{fig:coro_geo_variability}, this entails that the PC-USM-Net is not very effective in capturing the solution near the edge, where the solution is heavily affected by the geometric details of the domain not captured by the landmarks.
The use of a UC system is helpful in this regard by defining a model that receives as input, not the physical coordinates, but the reference ones.
These coordinates directly encode details of the geometry not captured by the landmarks.
In particular, the UC system makes the points belonging to the boundary of the various domains correspond to each other.
In this way, UC-USM-Nets are more effective than PC-USM-Nets in capturing the velocity field close to the boundary.

\begin{figure*}
    \centering
    \includegraphics[width = \textwidth]{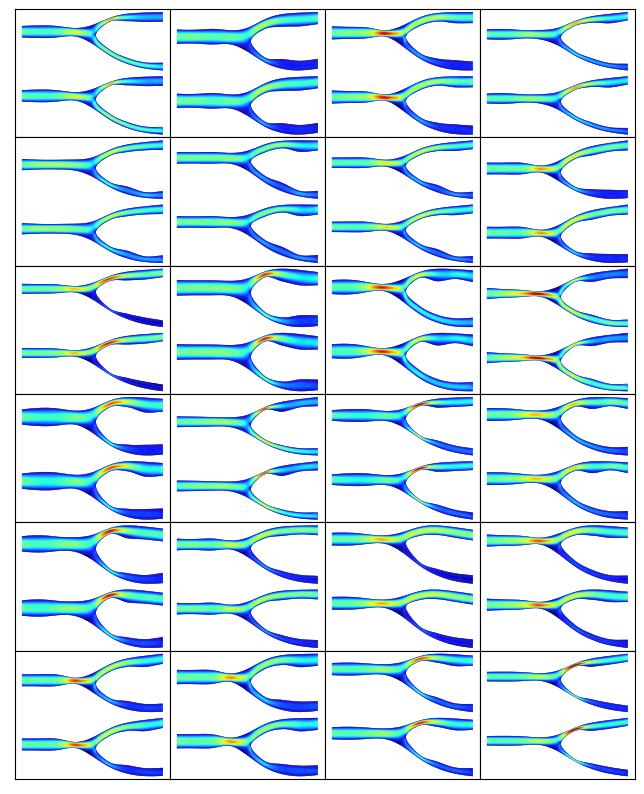}\\
    \includegraphics[]{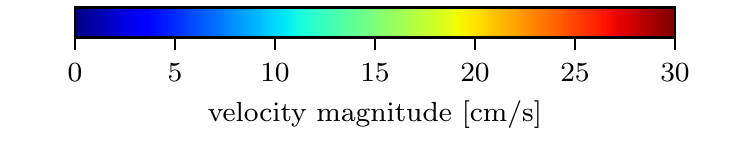}
    \caption{Test Case 2: comparison of the velocity magnitude field obtained with the FOM (top figure within each box) and with the UC-USM-Net (bottom figure within each box) in a subset of the test set.}
    \label{fig:coro_comparison_u}
\end{figure*}

\begin{figure*}
    \centering
    \includegraphics[width = \textwidth]{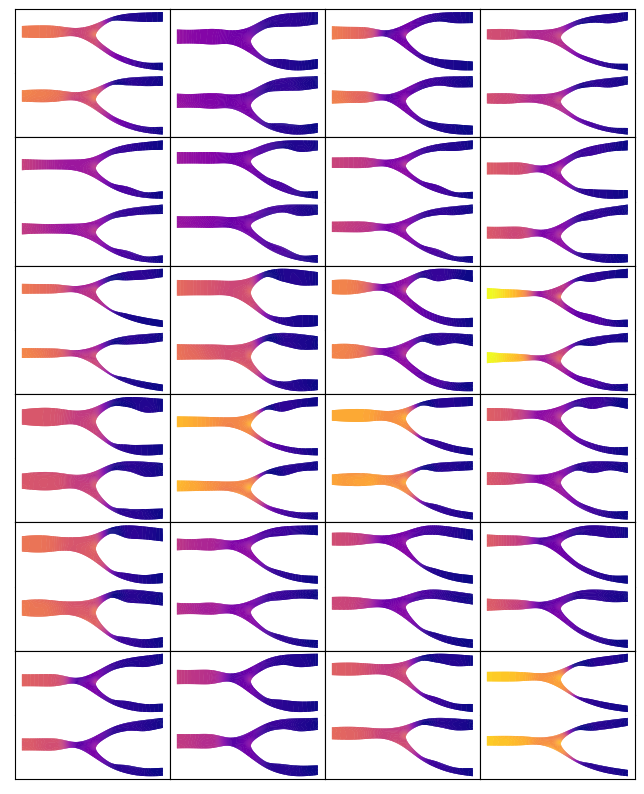} \\
    \includegraphics[]{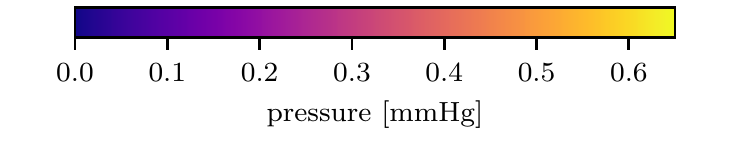}
    \caption{Test Case 2: comparison of the pressure field obtained with the FOM (top figure within each box) and with the UC-USM-Net (bottom figure within each box) in a subset of the test set.}
    \label{fig:coro_comparison_p}
\end{figure*}

In Figs.~\ref{fig:coro_comparison_u} and \ref{fig:coro_comparison_p} we show the velocity and pressure fields predicted by one of the trained UC-USM-Nets on a subset of the test dataset.

\section{Discussion} \label{sec:discussion}

We have introduced USM-Nets, a deep learning class of surrogate models capable of learning the solution manifold of a PDE universally with respect to physical parameters and geometry.

The ability of a surrogate model to capture the geometrical variability of the solution of a differential problem is a feature of great interest.
Indeed, many applications require considering the solution of a physical problem in different domains.
Biomedicine offers several examples in this regard since each patient presents a different geometry, and in many cases (as in hemodynamics) the geometry itself is the principal determinant of the solution.
Examples are given by the blood flow in an aneurysm, in a stenotic artery, or through an artificial valve.

Nevertheless, most of the reduced-order/surrogate models available in the literature consider a fixed domain accounting only for the variability of physical parameters \cite{benner2005dimension,antoulas2000survey,lassila2014model,benner2015survey,peherstorfer2015dynamic}. Few models rely on parametrized shape models that guarantee correspondence between points coming from different shapes, enabling the construction of projection-based models.
As a matter of fact, representing a solution manifold in variable geometries is an arduous task. There are two main difficulties in this regard: (1) how to encode the properties of the geometry at hand and (2) how to construct a discrete representation of the solution that is universal with respect to the shape of the domain.

Concerning point (1), USM-Nets only require the definition of a finite set of scalar quantities, called geometrical landmarks, that characterize the salient properties of the geometry at hand. Landmarks make USM-Nets an extremely flexible technique that can address a wide range of real-world applications.
There are different approaches to landmark definition, such as the one based on the statistical analysis of sampled geometries, like the first coefficients of proper orthogonal decomposition (POD) or the positions of control points.
However, approaches combining POD and ANN \cite{hesthaven2018non,carlberg2019recovering,dal2020data,o2021adaptive} might present difficulties in the database construction and limited generalization properties imposed by both the shape model, encoding the correspondence between points belonging to different geometries, and the truncation of the expansion.
Landmarks could be simply the coordinates of some points that characterize the geometry at hand. This case, specifically, is well suited for clinical applications.
Landmarks, such as the coordinates of a bifurcation, the position of an inlet, diameters, or areas can be processed directly from medical images without the need for segmentation and the generation of computational grids.
The great flexibility of USM-Nets lies in the fact that no structural requirements are imposed on the definition of landmarks.

Concerning point (2), a key feature of USM-Nets is their mesh-less nature, which frees them from a predetermined triangulation of the domain, overcoming the technical difficulties related to mesh element deformations.
The mesh-less nature of USM-Nets is achieved by their architectural design.
Unlike many existing surrogate modeling methods, that provide a map from the problem parameters to a set of degrees of freedom associated with a preconstructed parametrization of the solution trial manifold, the output of USM-Nets is the solution itself evaluated in a query point.
Indeed, by fixing a given parameter vector $\parVec$ and a given geometry $\Omega$, the USM-Net is a function from $\mathbb{R}^{\dimSpace}$ to $\mathbb{R}^{\dimSol}$, that is an approximation of the solution $\sol(\cdot; \parVec, \Omega)$.
Hence, instead of passing through a parametrization of the approximate solution, we make the ANN \textit{coincide} with the approximate solution.
A further advantage of this architectural design is that USM-Nets encode by construction the spatial correlation (that is, with respect to the input $\pt$) and do not need to learn it, thus achieving elevated accuracy levels even with lightweight NNs.

We have then presented an enhanced version of PC-USM-Nets, called UC-USM-Nets, based on a universal coordinate system.
Even if it is not straightforward in all practical cases to define a UC system (such as when the domain may vary in topology), the use of a UC system can improve the generalization accuracy of PC-USM-Nets, as shown by the numerical results.
A UC system acts at two levels.
Firstly, it allows us to partially compensate for the possible non-exhaustiveness of the geometrical landmarks in describing the geometry (see also point (1) of the discussion).
In Test Case 2, for example, in the setting with only six landmarks, we can have two geometries $\Omega^1$ and $\Omega^2$ that differ from each other, even though they have the same landmarks ($\geoMapping(\Omega^1) = \geoMapping(\Omega^2)$).
In boundary areas far from the landmarks, the PC-USM-Net might fail in satisfying the no-slip condition, while the UC-USM-Net, on the other hand, allows the solution to be more accurate because the UC system \textit{informs} the model of the position of the boundary in the geometry at hand.
Furthermore, regarding point (2), a UC system provides a more effective representation of the solution manifold.
In fact, in this case, the FCNN does not receive as input the coordinates $\pt \in \Omega$, but rather $\ptRef \in \DomainRef$, which are more informative of the \textit{role} that each point plays within the specific domain.

\section{Conclusions} \label{sec:conclusions}

We have proposed a novel technique (USM-Nets), based on ANNs, to build data-driven surrogate models that approximate the solution of differential equations while accounting for the dependence on both scalar physical parameters and the domain geometry.
Our method is non-intrusive as it does not require the knowledge of the FOM equations, but rather it is trained with samples of precomputed solution snapshots obtained for different parameters and geometries. 
It is also meshless since the USM-Net learns the map from point coordinates to the solution.
To characterize the geometrical features of the domain at hand, we consider a set of geometrical landmarks defined by the user.
Our method is highly flexible, as it does not pose specific requirements on the definition of these landmarks, making it suitable for practical applications and significantly easing its technological or clinical translation.

We have then presented an enhanced version of our surrogate modeling method, based on a UC system employed to pre-process the physical coordinates.
As shown by our numerical results, using this UC system enhances the generalization accuracy in some cases.

We have finally presented two test cases in fluid dynamics. The first is a lid-driven cavity problem with variable geometry and variable Reynolds number; the second one consists in predicting the steady-state pressure and velocity field within a coronary bifurcation, given the patient geometry.
In both test cases, despite the noticeable variability of the physical and/or geometrical parameters, USM-Nets were able to approximate the solution within an approximation of the order of $10^{-2}$, being trained, and a few hundreds of solution snapshots.


\bibliographystyle{asmems4}

\bibliography{references}


\end{document}